\documentclass[11pt,a4paper]{article}
\usepackage{amssymb}
\usepackage{amsthm}
\usepackage{amsmath}
\usepackage{amsbsy}
\usepackage[utf8]{inputenc}
\usepackage[a4paper, left=1.5cm, right=1.5cm, top=3cm, bottom=3cm]{geometry}
\usepackage[T1]{fontenc}
\usepackage{lmodern}
\usepackage{extarrows}
\usepackage{graphicx}
\usepackage{caption}
\usepackage{subcaption}
\usepackage[usenames,dvipsnames]{xcolor}
\usepackage[notcite,notref]{showkeys}

\makeatletter
\newtheorem*{rep@theorem}{\rep@title}
\newcommand{\newreptheorem}[2]{%
\newenvironment{rep#1}[1]{%
 \def\rep@title{#2 \ref{##1}}%
 \begin{rep@theorem}}%
 {\end{rep@theorem}}}
\makeatother

\newtheorem{theorem}{Theorem}[section]
\newreptheorem{theorem}{Theorem}
\newtheorem{definition}[theorem]{Definition}
\newtheorem{lemma}[theorem]{Lemma}
\newreptheorem{lemma}{Lemma}
\newtheorem{obs}[theorem]{Observation}
\theoremstyle{definition}

\DeclareMathOperator{\RE}{Re}

\DeclareMathOperator{\SPAN}{span}

\providecommand{\spann}[1]{\SPAN\left\{#1\right\}}
\providecommand{\boldU}{\boldsymbol{u}}

\providecommand{\zetKwadrat}{\mathbb{Z}^2}
\renewcommand{\d}{\partial}
\newcommand{\T}{\mathbb{T}}
\newcommand{\R}{\mathbb{R}}
\newcommand{\Z}{\mathbb{Z}}

\begin{document}
\begin{center}
{\Large \bf A construction of two different solutions to an elliptic system}

 \vskip 0.5cm
{\large
Jacek~Cyranka$^{*,\ddag}$, Piotr~Bogus\l aw~Mucha$^*$
}
\vskip 0.5cm
{\small$^*$ Institute of Applied Mathematics and Mechanics,
University of Warsaw}\\
{\small  Banacha 2, 02-097 Warszawa, Poland}
\vskip 0.5cm
{\small$^\ddag$ Department of Mathematics,
Rutgers, The State University of New Jersey}\\
{\small  110 Frelinghusen Rd, Piscataway, NJ  08854-8019, USA}

\vskip 0.5cm
{cyranka@mimuw.edu.pl, p.mucha@mimuw.edu.pl}
\vskip 0.5cm
\today
\end{center}

\providecommand{\matrixS}{{T}}
\providecommand{\seriesD}[2][l]{d_{#2}^{#1}}
\providecommand{\seriesA}[2][l]{a_{#2}^{#1}}
\providecommand{\seriesB}[2][l]{b_{#2}^{#1}}
\providecommand{\seriesAconst}{2^{2m}}
\providecommand{\operS}{{L_{\lambda}^l}}
\providecommand{\operSS}[2][\lambda]{{L_{#1}^{#2}}}

 {\bf Abstract.} The paper aims at constructing two different solutions to an elliptic system
 $$
u \cdot \nabla u + (-\Delta)^m u = \lambda F
 $$
 defined on the two dimensional torus.
 It can be viewed as an elliptic regularization of the stationary Burgers 2D system.
 A motivation to consider the above system comes from an examination of unusual propetries of the linear operator
 $\lambda \sin y \d_x w + (-\Delta)^{m} w$ arising from a linearization of the equation about the dominant part
 of $F$.
 We argue that the skew-symmetric part of the operator provides in some sense a smallness of norms of the linear operator inverse.
 Our analytical proof is valid for a particular force $F$ and for $\lambda > \lambda_0$, $m> m_0$ sufficiently large.
The main steps of the proof concern finite dimension approximation of the system and concentrate on
analysis of features of large matrices, which resembles standard numerical analysis.
Our analytical results are illustrated by numerical simulations.

 \vskip1cm

\paragraph{Keywords:} {nonlinear elliptic problem, 2D stationary Burgers equation, nonuniqueness, construction of two solutions, large matrices}
\paragraph{AMS classification:} {Primary: 35J60, 35A02. Secondary: 35Q99, 15B99.}

\section{Introduction}

Analysis of sets of solutions to elliptic systems/equations is of particular interest in the current research on partial differential equations (PDEs).
On the one hand, the question is challenging from the viewpoint of mathematical techniques. On the other hand,
the precise information about this set is crucial for understanding the dynamics of evolutionary problems
behind the elliptic one. In general, existing theory provides us with two answers: either there exists a single solution or
the system admits at least one solution.

Existing methods of the PDEs analysis provide only few examples for quite simple problems.
Starting with the classical example using the Mountain Pass Theorem for a semilinear elliptic equation \cite{E}.
Nonuniqueness example for the stationary Navier-Stokes equations \cite{G}, important geometric examples related to the
mean curvature problems \cite{BC}, or nonuniqueness for the one-dimensional viscous Burgers' equation \cite{BP12}, (in the case of 
evolutionary system we refer to
\cite{D}, \cite{AA}).
Derivation of asymptotic lower bound for the multiplicity of solutions for a semilinear problem can be found in \cite{DY}, 
\cite{LM}, and for a class of elliptic equations with jumping nonlinearities in \cite{MP}.
The work on numerical multiplicity proofs for systems/higher dimensional PDEs has also been an active topic of research.
There exist several
computer assisted proofs of existence of at least several solutions of certain parabolic PDEs. Let us stress that
contrary to our approach, all results obtained using a direct computer assistance are true essentially for some isolated parameter values or
a compact set of parameter values,
because all of the computations performed by a computer are finite.
Representative results include a proof of existence of four solutions to a semilinear
boundary value problem for particular choice of parameters \cite{BMP}, an outlook for the multiplicity of solutions for some
multidimensional PDEs is provided by a proof of existence of nonsymmetric solutions for a symmetric boundary value problem \cite{AK},
validated bifurcation diagrams constructed in \cite{BLV}, \cite{GL}, structure of the global attractor \cite{MMW}, 
numerical existence proofs for a fluid flow, and convection problems \cite{WN}, \cite{HNW}.

 \bigskip

 The subject of the present paper is the following elliptic system, which can be viewed as an elliptic regularization of the 
\emph{stationary Burgers system} \cite{B}, \cite{H}, \cite{Co} \emph{in 2D}
 \begin{equation}
   \label{eq:main}
 \boldsymbol{u} \cdot \nabla \boldsymbol{u} + (-\Delta)^m \boldsymbol{u} = \lambda \boldsymbol{F} \mbox{ \ \ on \ } \T^2.
 \end{equation}
Here $\boldsymbol{u}$ is sought as a vector function $\boldsymbol{u}:\T^2 \to \R^2$. The vector $\boldsymbol{F}$ is a particular external force,
and in this paper we define it as
\begin{equation}
\label{eq:F}
 \boldsymbol{F}(x,y)=\left(
 \begin{array}{c}
  \sin y \\
  \sin x
 \end{array}
\right).
\end{equation}
The magnitude of the external force is controlled by the parameter $\lambda$ and it is assumed to be greater than some positive number $\lambda_0$.
We shall note that the system has \emph{no a-priori estimate}. The issue of the existence of a solution to the system (\ref{eq:main}) is still open for a general form of
$\lambda \boldsymbol{F}$. To the best of our knowledge even the basic case of $m=1$ is unclear.


Let us discuss what motivated the presented research. Our numerical investigations of \eqref{eq:main}
revealed a solution possessing a curious structure: one of solution's Fourier modes being of $\lambda$ magnitude, and the remainder being bounded
uniformly with respect to $\lambda$. We further noticed that the natural symmetry embedded in this equation
implies the existence of a second solution, as the reflection by the symmetry of dominant part produces an essentially
different solution. Further on, to convince ourselves that this structure is in fact conserved for $\lambda$ large values after
a bifurcation, we performed a numerical bifurcation analysis, which showed that the graph of solution's norm is approximately
linear, and in fact there is a pitchfork bifurcation in the system for a particular $\lambda$.

We emphasize, that our situation is not as simple as it would be, when a symmetry embedded in the problem implies immediately existence of 
a different second solution. For sufficiently small $\lambda$'s solutions are symmetric, and for a large $\lambda$ the
symmetry is broken, which allows us to establish existence of two different solutions for large $\lambda$'s.
For small $\lambda$'s we claim only existence of a solution, as the two solutions from our main result merge into a single one.
The symmetry is seen as elementary, simply enough, we can exchange $x$ with $y$, and the first component with the second component (denoted $x\leftrightarrow y$ in the sequel).
Apparently, a stronger regularization effect is needed, than the one provided by the Laplacian operator. This is why 
we state our main result (Theorem~\ref{thm:main}) for $m$ sufficiently large. Our analytical results are supported by a numerical bifurcation analysis
(Section~\ref{sec:numBifAna}).

The main tool of our technique is to exploit unusual features of a linearization of the system.
Let $\|w\|_{l^\infty}$ denote the supremum norm of elements of the Fourier series $w$. Apparently,
for the solutions to the following scalar problem
\begin{equation}
  \label{eq:scalarProblem}
 \lambda \sin y \d_x w + (-\Delta)^{m} w = \lambda \sin x \mbox{ on } \T^2,
\end{equation}
we obtain
 \begin{equation}
  \|w\|_{l^\infty} \lesssim 1,
 \end{equation}
in other words, this quantity is free from $\lambda$ dependence (for large $\lambda$), although
other norms are growing with $\lambda$.
We see an interplay between the growth of the right-hand side and an increase of influence of the term $\lambda \sin \partial_x$, which represents
(in some sense) a rotation effects. In particular, it causes that the amplitudes of modes to be uniformly bounded. Such effect can be 
compared with general phenomenon of hypocoercivity explained in \cite{Vil09}. We shall note, however, we do not apply the general theory
for operators of type $A^*A+B$, since we want to avoid considerations in Hilbertian spaces. We work instead in the $l^1,l^\infty$
framework, which is the most optimal for our analysis. The technique is elementary, in order to obtain a constructive bounds for
linear operators inverses we perform a \emph{large matrices analysis}.
The features of the linear operator are first found for its finite dimensional truncation -- a Galerkin approximation, then
the properties of the full infinite dimensional operator are obtained using a limit passage.
The key result concerning (\ref{eq:scalarProblem}) is described by Theorem \ref{lem:invMatrixEstimates} and its proof is the main part of this paper.

Our analysis of the system \eqref{eq:main} allows to prove the following theorem being the main result of the present paper.

\smallskip

\begin{theorem} 
\label{thm:main}
Let $m>9/2$ and $\lambda > \lambda_0$ be sufficiently large. Then there exists at least two solutions to
the system \eqref{eq:main} with $F$ given by \eqref{eq:F} such that
\begin{equation}\label{i5}
  \mathfrak{u}_1 =
 \lambda \left(
 \begin{array}{c}
  \sin y \\
  0
 \end{array}
\right)
+ L^1(x,y) + R^1(x,y),
\end{equation}
and
\begin{equation}\label{i6}
 \mathfrak{u}_2 =
 \lambda \left(
 \begin{array}{c}
  0 \\
  \sin x
 \end{array}
\right)
+ L^2(x,y) + R^2(x,y),
\end{equation}
where $L^1,L^2$ are solutions to the linearization and they are of order $\lambda^{2/m}$ in $l^\infty$ norm and
$R^1,R^2$ are of order $\lambda^{-\alpha}$ with $\alpha >0$.
\end{theorem}

\smallskip

The proof of Theorem~\ref{thm:main} is based on a subtle analysis of the system (\ref{eq:scalarProblem}). We impose the form of solutions and then we construct them via approximation on
finite dimensional subspaces. The natural symmetry $x \leftrightarrow y$ implies that we obtain at least two different solutions, provided $\lambda$ is sufficiently large.

Indeed, the properties of the system \eqref{eq:scalarProblem} established in Theorem~\ref{thm:main} are the main impact of the present paper.
We are ensured that this type of properties will allow to study precise dynamics of systems
 with the transport term $u\cdot \nabla u$. The most natural example is the Navier-Stokes equations. However, the current box of tools is not sufficient to attack this problem.
 We present here a brand new technique to study quasilinear elliptic systems. Hence one can look at the system \eqref{eq:main} as a toy model
 for which we demonstrate our new method.

We are highly convinced that the explicit bounds for norms of tridiagonal differential operators obtained in this work,
which are independent of the dimension, can be applied to study other problems, including bounds for solutions of some
linear PDEs, computer assisted proofs for nonlinear PDEs, numerical analysis of discretizations of certain PDEs, 
and slow-fast systems. There are existing research efforts in understanding structure of the tridiagonal operators arising in PDEs, see e.g. \cite{BDL}.
Let us also note that methods based on Fourier series may be applied for systems in pipe-like domains. An example is \cite{M}, where analysis of the Oseen operator gave very precise space asympotics of solutions in front and behind an obstacle.

\smallskip

Literature concerning the issue of existence of solutions to the stationary Burgers equations is not rich. Most of the results concern
only the mono-dimensional case model \cite{BGS01,BP12,BSG11}. It motivates us to perform 
 numerical analysis for various cases of the model. We observe that for the system
\eqref{eq:main} there exist a threshold value $\hat{m}>1$, such that the main result is valid for all $m>\hat{m}$ --
the two distinct solutions can be still constructed.
For the case $m<\hat{m}$, especially, in case of the stationary forced 2D Burgers equations $(m=1)$ the global
picture is significantly different, and for $m=1$ certainly the two solutions cannot anymore be isolated as in
the other cases. 

\smallskip

We note that since finishing of the first version of this paper the first author significantly improved the
  bounds for the norm of the inverse tridiagonal operators \cite{CL}. The motivation has been to develop a
  validated numerical scheme for forward integration of a class of parabolic PDEs. We are now convinced that a proof along
  the lines presented in this paper is possible for the stationary Burgers system with smaller exponents $m$
  defining the linear operator in \eqref{eq:scalarProblem}.
  Moreover, we are convinced that a similar proof is also possible for other problems,
  including the stationary 2D viscous Navier-Stokes equations. We will investigate this possibilities in future research.

The paper is organized as follows. We present in Section~\ref{sec:preliminaries} the subject of this paper 
written in coordinates, in Section~\ref{sec:numBifAna} bifurcation diagrams, and a brief technical explanation. In 
Section~\ref{sec:symmetries}, the relevant symmetries of the problem, which are crucial in our analysis. 
In Section~\ref{sec:linearized}, the matrix form of the linearized operator, along with 
some important inverse operators bounds. In Section~\ref{sec:fixedPoint}, a-priori bounds for the
solutions of finite dimensional truncations, and in Section~\ref{sec:proof}, an existence argument
for the infinite dimensional system. Finally, in Section~\ref{sec:technicalLemmas}, some technical lemmas
necessary to prove crucial inverse operators bounds from Section~\ref{sec:linearized}.

\subsubsection*{Acknowledgments}

The presented work has been done while the first author held a post-doctoral
position at Warsaw Center of Mathematics and Computer Science, his
research has been partly supported by Polish National Science Centre grant 2011/03B/ST1/04780.
The second author (PBM) has been partly supported by National  Science  Centre  grant
2014/14/M/ST1/00108 (Harmonia).

\section{Preliminaries}
\label{sec:preliminaries}

We start our analysis with the preparation of our system
\begin{equation}
  \label{eq:mainSystem}
 \begin{array}{l}
  u \cdot  \nabla u^1 + (-\Delta)^m u^1 = \lambda \sin y, \\[6pt]
  u \cdot\nabla u^2 + (-\Delta)^m u^2 = \lambda \sin x.
 \end{array}
\end{equation}
We  fix the notation
\begin{equation}
  \label{ansatz}
 u=\bar v + V, \mbox{ \ \ \ where \ \ \ } \bar v = \lambda \left(
 \begin{array}{c}
  \sin y \\
  0
 \end{array}
\right) .
\end{equation}
We focus just on construction of solution (\ref{i5}), the symmetry will imply existence of the second one -- see Section 4 (Definition~\ref{defsymab}).
The above relations restate the system \eqref{eq:mainSystem} as follows
\begin{equation}\label{p1}
 \begin{array}{l}
  \lambda \sin y \d_x V^1 + (-\Delta)^m V^1 = - \lambda \cos y \, V^2 - V \cdot \nabla V^1, \\[6pt]
  \lambda \sin y \d_x V^2 + (-\Delta)^m V^2 = \lambda \sin x - V \cdot \nabla V^2.
 \end{array}
\end{equation}
Observe that the term $\lambda\sin{y}$ is not present in \eqref{p1}, as it disappears due to the ansatz \eqref{ansatz}.
In order to split the solution into two parts, the first with small amplitudes and the second with higher ones. We introduce a linearization of (\ref{p1})
\begin{equation}
  \label{eq:linearization}
 \begin{array}{l}
  \lambda \sin y \d_x A + (-\Delta)^m A = - \lambda \cos y \, B, \\ [6pt]
  \lambda \sin y \d_x B + (-\Delta)^m B = \lambda \sin x,
 \end{array}
\end{equation}
and define $V$ as the following pair
\begin{equation}
 V=\left(
  \begin{array}{c}
   A \\
   B
  \end{array}
  \right)+
  \left(\begin{array}{c}
   a \\
   b
  \end{array}
  \right).
\end{equation}
Vector $(A,B)^T$ defines $L^1$ appearing in (\ref{i5}).
  This step of prescription of constructing solutions to (\ref{eq:main}) is important, since (\ref{eq:linearization}) implies a constraint on $A$ and $B$.
  This relation turns out to be satisfied also by $a$ and $b$.
By differentiating $(\ref{eq:linearization})_1$ with respect to $x$ and $(\ref{eq:linearization})_2$ with respect to $y$, the system \eqref{eq:linearization} takes the form
\begin{equation}
  \label{eq:ABd}
 \begin{array}{l}
  \lambda \sin y \d_x A_x + (-\Delta)^m A_x = - \lambda \cos y \, B_x, \\ [6pt]
  \lambda \sin y \d_x B_y + (-\Delta)^m B_y = - \lambda \cos y \, B_x,
 \end{array}
\end{equation}
So we obtain 
\begin{equation}
 \lambda \sin y \partial_x(A_x -B_y) + (-\Delta)^m(A_x-B_y)=0.
\end{equation}
Testing it by $(A_x-B_y)$ we get
\begin{equation}
 \int_{\T^2} | \nabla (-\Delta)^{m/2-1}(A_x-B_y)|^2 dxdy=0, \mbox{ \ and of course \ } \int_{\T^2} (A_x-B_y) dxdy =0.
\end{equation}
Hence we get the desired constraint 
\begin{equation}\label{AB-con}
 A_x = B_y.
\end{equation}

Returning to $V$  we find equations for $a$ and $b$
\begin{equation}\label{p2}
 \begin{array}{l}
  \lambda \sin y \d_x a + (-\Delta)^m a = - \lambda \cos y \,  b - \left(
  \begin{array}{c}
   a + A \\
   b + B
  \end{array}
  \right) \cdot \nabla (a+A), \\
  \lambda \sin y \d_x b + (-\Delta)^m b =  - \left(
  \begin{array}{c}
   a + A \\
   b + B
  \end{array}
  \right) \cdot \nabla (b+B).
 \end{array}
\end{equation}

Here again one can check constraint (\ref{AB-con}) for $a$ and $b$. Taking suitable differentiation  of system (\ref{p2}) we find
\begin{equation}\label{p2a}
 \begin{array}{l}
  \lambda \sin y \d_x a_x + (-\Delta)^m a_x = - \lambda \cos y \,  b_x - \left(
  \begin{array}{c}
   a_x + A_x \\
   b_x + B_x
  \end{array}
  \right) \cdot \nabla (a+A)
  - \left(
  \begin{array}{c}
   a + A \\
   b + B
  \end{array}
  \right) \cdot \nabla (a_x+A_x), \\
  \lambda \sin y \d_x b_y + (-\Delta)^m b_y = - \lambda \cos y\,b_x - \left(
  \begin{array}{c}
   a_y + A_y \\
   b_y + B_y
  \end{array}
  \right) \cdot \nabla (b+B)
  - \left(
  \begin{array}{c}
   a + A \\
   b + B
  \end{array}
  \right) \cdot \nabla (b_y+B_y).
 \end{array}
\end{equation}
So then we find, keeping in mind (\ref{AB-con})
\begin{equation}\label{p2b}
 \lambda \sin y \partial_x(a_x-b_y) +(-\Delta)^m (a_x-b_y)= - \left[ (a_x+A_x+B_y+b_y)(a_x-b_y)\right] - \left(
  \begin{array}{c}
   a + A \\
   b + B
  \end{array}
  \right) \cdot \nabla (a_x-b_y) .
\end{equation}

Observe that as the rhs of (\ref{p2b}) would be zero than we find the 
desired  constraint 
\begin{equation}\label{ab-con}
 a_x=b_y.
\end{equation}
This relation will be guaranteed by the construction presented at the beginning of Section 6.
In few words, the construction is performed via an iteration scheme, so vanishing of the rhs will be guaranteed by the previous step, see (\ref{f13a}).

Looking at the above problems we see that the analysis depends on the properties of the following operator
 \begin{equation}
 L_\lambda(w)=\lambda \sin y \d_x w +(-\Delta)^m w.
\end{equation}

The key element of the proof of Theorem~\ref{thm:main} is a result concerning norm estimates for the $L_\lambda$ inverse operator. The precise statement of the result we find in Section~\ref{sec:invOperatorBounds}, it is Theorem~\ref{thm:main}.

\paragraph{Notation} In bold we denote complex coefficients, e.g. $\boldsymbol{a}_k=\left(a_k^1,a_k^2\right)\in\mathbb{C}^2$,
where $a_k^1$, and $a_k^2$ denote the first, and the second component of $\boldsymbol{a}_k$ respectively.
Let $k,k_1,k_2\in\mathbb{Z}^2$ denote pairs of integers. By $k^1,k^1_1,k^1_2$, and $k^2,k^2_1,k^2_2$ we denote the first, and the second components respectively.

\paragraph{}We rewrite the problem \eqref{eq:main} using Fourier's coordinates, being the most natural way to consider problems on a torus.
\begin{subequations}
  \label{eq:toyInFourier}
  \begin{gather}
    \boldU(x)=\sum_{k\in\mathbb{Z}^2}{\boldsymbol{a}_ke^{ik\cdot(x,y)}},\quad \boldsymbol{a}_k=(a_{k}^1,a_{k}^2)\in\mathbb{C}^2,\\
    \boldsymbol{F}(x)=\sum_{k\in\mathbb{Z}^2}{\boldsymbol{F}_ke^{ik\cdot(x,y)}},\quad \boldsymbol{F}_k=(F_{k}^1,F_{k}^2)\in\mathbb{C}^2.\\
    \sum_{\substack{k_1+k_2=k\\k_1,\,k_2\in\mathbb{Z}^2}}{a_{k_1}^1ik_{2}^1a_{k_2}^j}+\sum_{\substack{k_1+k_2=k\\k_1,\,k_2\in\mathbb{Z}^2}}{a_{k_1}^2ik_{2}^2a_{k_2}^j}+\left((k^1)^{2m}+(k^2)^{2m}\right)a_{k}^j-\lambda F_{k}^j=G(\boldsymbol{a},\lambda)_k^j=0,\quad j=1,2,\quad k\in\zetKwadrat\label{eq:infDim}.
  \end{gather}
The operator $(-\Delta)^m$ is diagonal in Fourier's basis,
having $\left((k^1)^2+(k^2)^2\right)^m$ as the eigenvalues. In order to simplify the arguments in the remainder of the paper as the operator
$(-\Delta)^m$ we will consider an operator having $(k^1)^{2m}+(k^2)^{2m}$ as the eigenvalues.
Of course, all of the presented arguments are also valid for the original  case, as $\left((k^1)^2+(k^2)^2\right)^m$ clearly bounds
$(k^1)^{2m}+(k^2)^{2m}$ from above.

For the particular choice of the external forcing, $\boldsymbol{F}$ is given by
\begin{equation}
  \tag{1d}
  F^1_{(0,1)}=F^2_{(1,0)} = -\frac{i}{2},\quad F^1_{(0,-1)}=F^2_{(-1,0)} = \frac{i}{2},\quad F^j_{k}=0\text{ for all the other cases.}
\end{equation}
\end{subequations}

\begin{definition}
  In the space of complex sequences $\{a_k\}_{k\in\mathbb{Z}^d}$, we will say that the sequence $\{a_k\}$ satisfies the reality
condition iff
\begin{equation}
  \label{eq:realityCondition}
  a_{k}=\overline{a_{-k}},\quad k\in\mathbb{Z}^d.
\end{equation}
\end{definition}

In the considered problem we impose odd periodic boundary conditions, i.e.
\begin{equation}\label{odd}
  \left.\begin{array}{lll}
      u^j(x,y)&=-u^j(-x,-y)&\ \\u^j(x,y)&=u^j(x+2\pi,y)&=u^j(x,y+2\pi)\end{array}\right.j=1,2,\quad x,y\in\mathbb{R},
\end{equation}
which on the level of the Fourier series means that we restrict the basis to odd functions, or
equivalently the coefficients are purely imaginary numbers satisfying
\begin{equation}
  \label{eq:fourierSines}
  \RE(a_{k}^j)=0,\ a_k^j = -a_{-k}^j\quad j=1,2,\ k\in\zetKwadrat.
\end{equation}
It is immediately verified that the space of coefficients satisfying \eqref{eq:fourierSines} is
invariant under the equation \eqref{eq:toyInFourier}, and we skip the formal calculations.
Observe that \eqref{eq:fourierSines} together with the reality condition implies automatically the following 'zero mass' constraint
\begin{equation}
  a^j_0=0,\quad j=1,2.
\end{equation}
Immediately, also we recognize that symmetry (\ref{odd}) implies that our solutions will be constructed as series in sinus only.

From now on we are going to consider the following finite dimensional approximation of the system \eqref{eq:toyInFourier}
\begin{definition}
Let $N>0$. We call the $N$-th Galerkin approximation of \eqref{eq:toyInFourier} the following system
  \begin{equation}
    \label{eq:projection}
    \tag{1P}
    \sum_{\substack{k_1+k_2=k\\|k_1|,\,|k_2|\leq N}}{a_{k_1}^1ik_{2}^1a_{k_2}^j}+\sum_{\substack{k_1+k_2=k\\|k_1|,\,|k_2|\leq N}}{a_{k_1}^2ik_{2}^2a_{k_2}^j}+\left((k^1)^{2m}+(k^2)^{2m}\right)a_{k}^j-\lambda F_{k}^j=G_N(\boldsymbol{a},\lambda)_k^j=0,\quad j=1,2,\  |k|\leq N.
  \end{equation}
\end{definition}

\begin{definition}
\label{banach}
 We introduce the following Banach spaces for sequences $\{a_k\}$ equipped with the following norms
 \begin{equation}
  \begin{array}{ll}
   \displaystyle\|\{a_k\}\|_{l^\infty}=\sup_{k\in \Z^2} |a_k|,  & \displaystyle\|\{a_k\}\|_{l^\infty_p}=\sup_{k\in \Z^2} (|k_1|+|k_2|)^p|a_k|, \\[10pt]
   \displaystyle\|\{a_k\}\|_{l^1}=\sum_{k\in \Z^2} |a_k|, & \displaystyle\|\{a_k\}\|_{l^1_1}=\sum_{k\in \Z^2} (|k_1|+|k_2|)|a_k|.
  \end{array}
 \end{equation}
\end{definition}

The norm used for multi-indices is taken to be the $\infty$ norm, i.e.
\[
  |k| := \max\left\{|k^1|, |k^2|\right\}.
\]

\begin{definition}
  Let us define the following space
  \begin{equation*}
    H = H_{(N)} = \left\{\{a_k\}\in\mathbb{C}^{(2N+1)^2-1}\colon a_k=\overline{a_{-k}},\ \RE(a_{k})=0\text{, for }0<|k|\leq N\right\}.
  \end{equation*}
  we are going to look for solutions $\boldsymbol{a}$ of \eqref{eq:projection}, such that $\boldsymbol{a}\in H\times H$.
  In the sequel, whenever $H$ appears, $N$ will either be fixed or clear from context.
\end{definition}

\section{Numerical bifurcation analysis}
\label{sec:numBifAna}

\begin{figure}[htbp]
  \centering
          \includegraphics[width=0.6\textwidth]{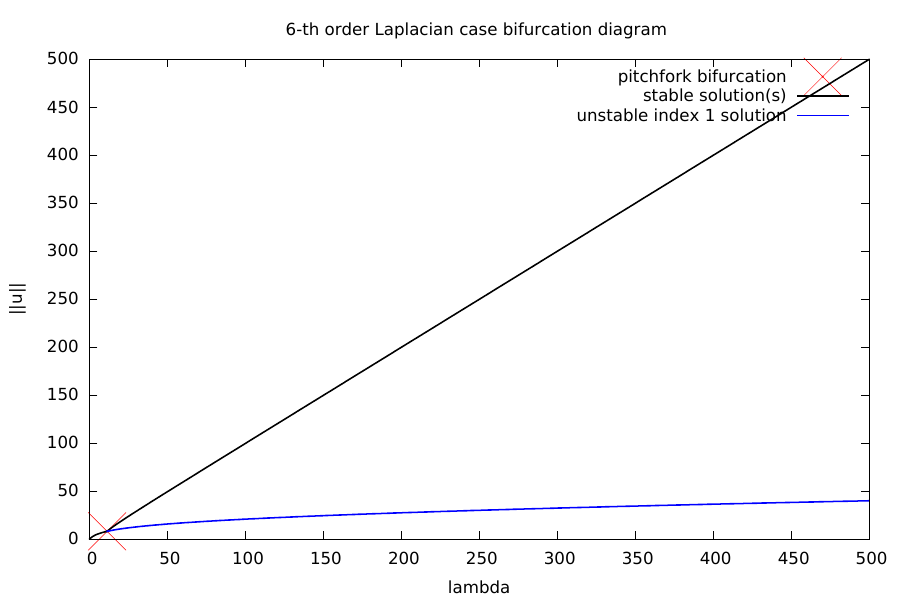}
          \caption{Bifurcation diagram for $u\cdot\nabla u+(-\Delta)^6u=\lambda F$. Clearly, the stable (black) solutions curve is
          almost $\|u(\lambda)\|=\lambda$, as there is only one mode of $\lambda$ magnitude. Single pitchfork bifurcation was
          detected at $\lambda_0=8.0629$.}          
  \label{fig:bif}
\end{figure}

We analyze the bifurcation structure of the problem \eqref{eq:projection}, and we present the results on Figure~\ref{fig:bif}.
Starting from the zero solution at $\lambda=0$ we follow the branch of solutions.
We detected a pitchfork bifurcation at a value of $\lambda$, which depends on the parameter $m$ appearing in \eqref{eq:projection}.
From the point of the pitchfork bifurcation we follow both the stable (one of two) and unstable branch (it is unique).

For a given $\lambda$ we solve for $a(\lambda)$ such that $G_N(a(\lambda),\lambda)=0$. We
implemented a path following procedure in order to track $a(\lambda)$.
To make any path following procedure work the partial derivative 
$\frac{\partial a(\lambda)}{\partial \lambda}$ is required, as bifurcation points
are detected by monitoring for its eigenvalues crossing zero. We implemented our path following
procedure on the top of the existing C++ software \cite{C} in which the partial derivative is calculated by
means of automatic diffrentiation and fast Fourier transforms, refer to \cite{C} for details.

\begin{figure}[htbp]
  \centering
        \begin{subfigure}[b]{0.4\textwidth}
          \includegraphics[width=\textwidth]{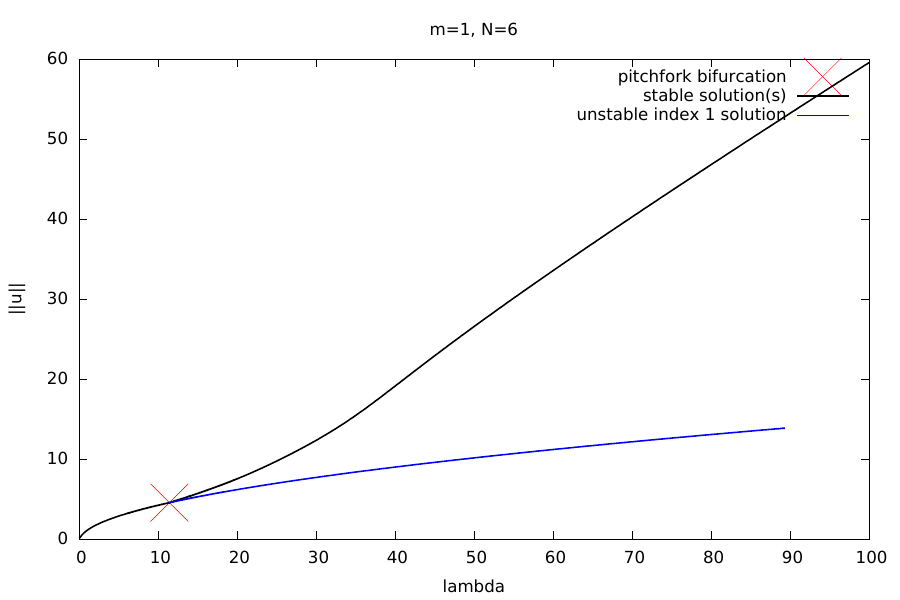}
        \end{subfigure}
        \begin{subfigure}[b]{0.4\textwidth}
          \includegraphics[width=\textwidth]{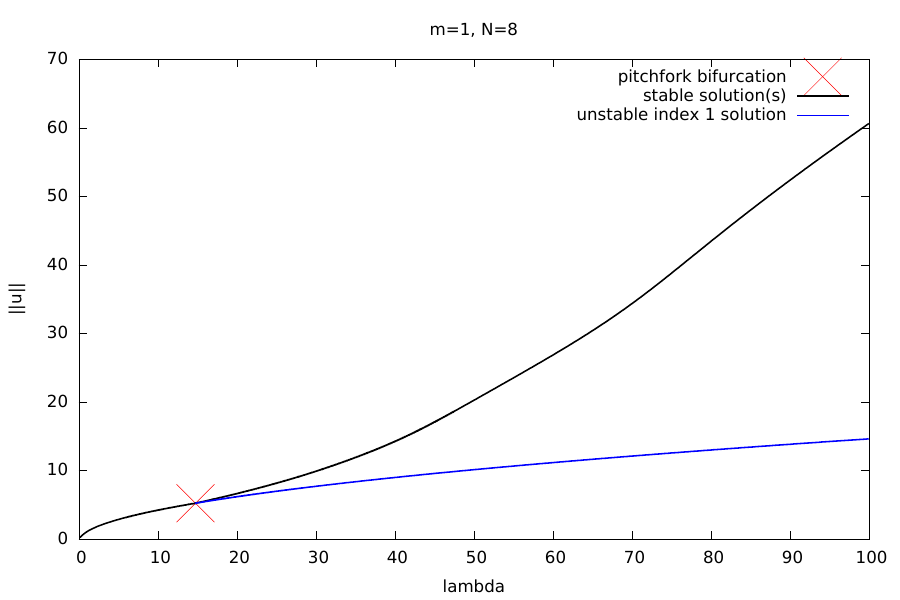}
        \end{subfigure}
        \begin{subfigure}[b]{0.4\textwidth}
          \includegraphics[width=\textwidth]{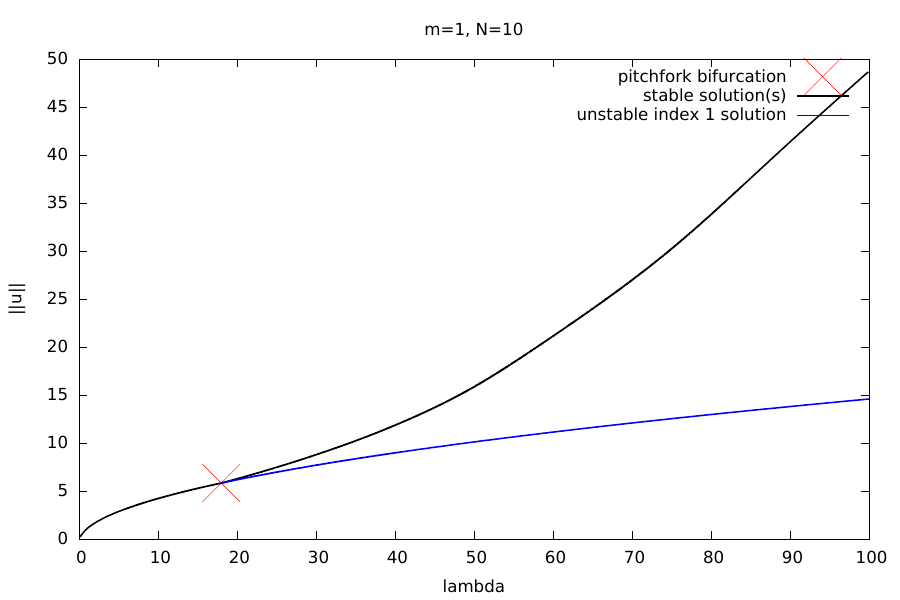}
        \end{subfigure}
        \begin{subfigure}[b]{0.4\textwidth}
          \includegraphics[width=\textwidth]{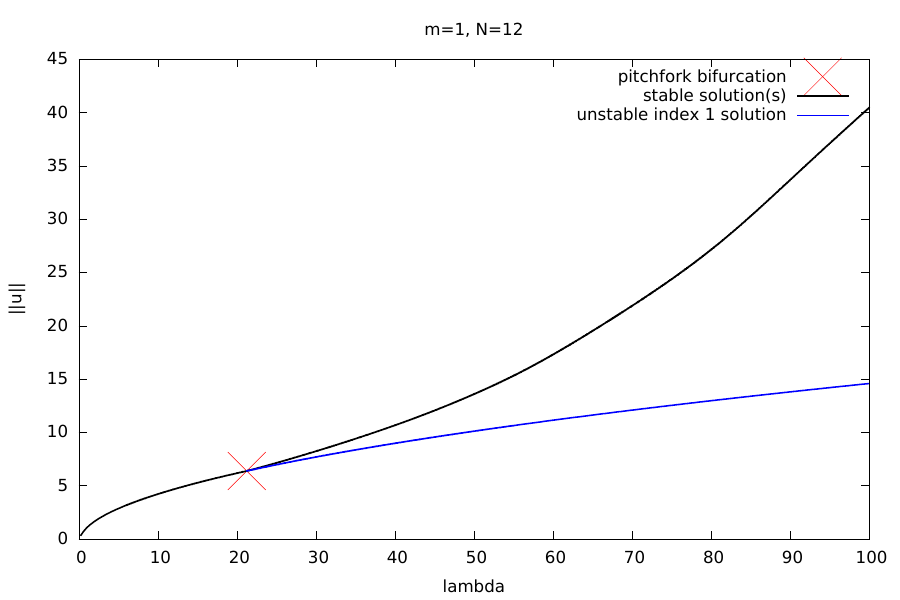}
        \end{subfigure}
  \caption{Bifurcation diagram for $u\cdot\nabla u+(-\Delta)u=\lambda F$. Each diagram was computed with 
different approximation dimension $N$ (given in the title). For this problem, the linear structure is not anymore evident, 
the question of the existence of two distinct solutions is left open in this case. The bifurcation point depends heavily
on the dimension.}
  \label{fig:bifLaplacian}
\end{figure}

We computed bifurcation diagrams for two specific cases
\begin{itemize}
  \item Figure~\ref{fig:bifLaplacian}, $m=1$ is fixed, and the truncation dimension $N$ is varied.
This case is excluded from our theory. 
\item Figure~\ref{fig:bif}, $m=6$ is fixed, and the truncation dimension $N$ is varied.
  This case is excluded from our theory.
\end{itemize}

There are some apparent differences between those two cases.
In Figure~\ref{fig:bifLaplacian} and \ref{fig:bif} in blue we marked the unstable branch of index $1$, and in black the 
stable solution(s) -- this branch represents in fact two solutions having the same norm related with a symmetry. 
The symmetry is denoted by $S$ in Section~\ref{sec:symmetries}.
Apparently, the considered pitchfork bifurcation is the point where the symmetry $S$ breaks.
Let us relate the presented diagrams with our theoretical results presented in the sequel. We prove that on the stable branch in
Figure~\ref{fig:bif} there are two distinct solutions, and this branch is approximately linear with respect to $\lambda$ for 
sufficiently large $\lambda$.

The diagrams were generated using the approximation with $N=8$, corresponding to $17^2/2$ degrees of freedom.

We approximate the solution using a fixed number of Fouriers' functions.
On Figure~\ref{fig:bifLaplacian} we present a few bifurcation diagrams obtained using Fouriers' approximation
with varying approximation dimensions (limited by our computational resources).
To construct the diagrams, we start from the zero solution at $\lambda=0$, the branch of solutions ($u(\lambda)$) is
followed until a bifurcating solution is found. 
In case a bifurcating solution is found, both of the branches: the original, and the new bifurcating
branch are followed.

Observe that those diagrams significantly differ. 
For instance the value of $\lambda$ for which the numerical pitchfork bifurcation occurs is proportional
to the approximation dimension, we mean that $\lambda$ is significantly larger, when a larger approximation
dimension is used. This leads us to the conjecture that the
apparent bifurcation is only a numerical artifact. It appears that in the case of stationary forced 2D Burgers 
equations $(m=1)$ the dynamics is either not finite dimensional, or the dimension of the attractor is really high. 
This is in contrary to the cases
included in our theory (e.g. $m=6$), in which the dynamics is essentially finite dimensional (the bifurcation
diagrams computed using different approximation dimensions does not
differ much). One possible explanation is that in case $m=1$ the
Laplacian operator does not provide strong enough smoothing effect compared to the higher order elliptic operators.

\providecommand{\symSpace}{H^\prime}
\providecommand{\symS}{S}
\providecommand{\N}{\mathbb{N}}

\section{Definition of a subspace of symmetric solutions}
\label{sec:symmetries}
In this section we define the symmetry exhibited by the studied problem, and which we will use in the sequel.
We make a standing assumption that the external forces that we consider are also symmetric.
Later on it will became evident that the second solution is obtained through the reflection by the  symmetry $\symS^{x\leftrightarrow y}$
(Definition~\ref{defsymab}).

Recall our working space of sequences of complex Fourier modes satisfying the reality condition
\begin{equation*}
    H = H_{(N)} = \left\{\{a_k\}\in\mathbb{C}^{(2N)^2+1}\colon\ \RE(a_{k})=0,\ a_k = -a_{-k}\right\}.
\end{equation*}

Instead of working directly with the space $H$, we will work with
the following product space of sequences of complex Fourier modes satisfying certain symmetry
exhibited by the solutions of the system \eqref{eq:toyInFourier}.

\begin{definition}
  \label{defsyms}
  Let $\symSpace$ be the following space
  \[
    \symSpace_{(N)} = \left\{\mathbf{a}\in H_{(N)}\times H_{(N)}\colon\,\mathbf{a}\text{ satisfies }\symS\,\mathbf{a} = \mathbf{a}\, \right\}.
  \]
\providecommand{\even}{\mathbb{Z}_{even}}\providecommand{\odd}{\mathbb{Z}_{odd}}
\emph{The symmetry $\symS\colon H\times H\to H\times H$} is the following symmetry. We define the symmetry directly
on the level of Fourier modes
\begin{equation}
\label{eq:symmetry3}
\begin{array}{lll}
  \symS\left(a^1_{(k^1,k^2)}, a^2_{(k^1,k^2)}\right) = \left(-a^1_{(-k^1,k^2)}, a^2_{(-k^1,k^2)}\right)&\text{ for }k^1,\,k^2\in\even\text{, or }k^1,\,k^2\in\odd,\\
  \symS\left(a^1_{(k^1,k^2)}, a^2_{(k^1,k^2)}\right) = \left(a^1_{(-k^1,k^2)},-a^2_{(-k^1,k^2)}\right)
  &\text{ for }k^1\in\even,k^2\in\odd\text{, or }k^1\in\odd,k^2\in\even.
\end{array}
\end{equation}


In Lemma~\ref{lemsym} we show that the nonlinearity appearing in the system \eqref{eq:toyInFourier} is symmetric with respect to $\symS$,
i.e. $\boldsymbol{G}(\boldsymbol{a})=0\iff \boldsymbol{G}(\symS\boldsymbol{a})=0$ as long as $\boldsymbol{F}$ is symmetric.

\end{definition}

\begin{obs}
Using the isomorphism of the space of sequences of Fourier modes with the space of functions
spannded by the trigonometric basis, space $H_{(N)}$ is isomorphic to a space of functions spanned by sines, i.e.
\begin{equation}
  \label{eqHprime}
 H_{(N)} \sim \left\{ \sum_{\substack{-N\leq l\leq N\\0\leq k\leq N}} v_{lk} \sin(lx+ky)\right\}.
\end{equation}
\end{obs}

\providecommand{\Nsymb}{Nonl}
\begin{lemma}
  \label{lemsym}
  Let $\boldsymbol{\Nsymb}$ be the nonlinear part of \eqref{eq:infDim} modulo the imaginary unit factor. $\boldsymbol{\Nsymb}$ satisfies
  \begin{equation*}
    \boldsymbol{\Nsymb}(\symS\boldsymbol{a})=\symS\,\boldsymbol{\Nsymb}(\boldsymbol{a}).
  \end{equation*}
\end{lemma}
\paragraph{Proof}
\[
  \boldsymbol{\Nsymb}(\boldsymbol{a}) = \left(\Nsymb^1(\boldsymbol{a}), \Nsymb^2(\boldsymbol{a})\right).
\]

Below, we check that the first component $\Nsymb^1$ satisfies the symmetry, by the same arguments the symmetry of the second component $\Nsymb^2$ follows.
To verify the claim let us consider two subcases
\paragraph{Case 1}
$k=(k^1,k^2)$, $k^1,\,k^2$ even or $k^1,\,k^2$ odd.
\begin{multline*}
  \Nsymb^1(\boldsymbol{a})_{(-k^1,k^2)}=\sum_{k_1+k_2=k}{a^1_{(-k_1^1,k_1^2)}(-k_2^1)a^1_{(-k^1_2,k^2_2)}}+\sum_{k_1+k_2=k}{a^2_{(-k^1_1,k_1^2)}k^2_2a^1_{(-k^1_2,k^2_2)}}=\\
  -\sum_{k_1+k_2=k}{(S\boldsymbol{a})^1_{(k_1^1,k_1^2)}k_2^1(S\boldsymbol{a})^1_{(k^1_2,k^2_2)}}-\sum_{k_1+k_2=k}{(S\boldsymbol{a})^2_{(k^1_1,k_1^2)}k^2_2(S\boldsymbol{a})^1_{(k^1_2,k^2_2)}}=\\
  -\Nsymb^1(\symS\boldsymbol{a})_{(k^1,k^2)} = \Nsymb^1(\symS\boldsymbol{a})_{(-k^1,k^2)}.
\end{multline*}
If we consider indices $k_1,\,k_2$ such that $k=k_1+k_2$, it holds that either $k^1_j,\,k^2_j$ are even (odd) ($k^1,\,k^2$ even case) or
one of $k^1_j,\,k^2_j$ is even and the second one is odd ($k^1,\,k^2$ odd case), $j=1,2$. This implies the second equality above, where in the first term
the symmetry generates either none or two minuses, as both of the modes come from the same component, hence, 
the only minus appears in front of the index $-k_2^1$. Whereas in the second term there is single minus generated, as the modes come from
different components, this is seen clearly from \eqref{eq:symmetry3}.
\paragraph{Case 2}
$k=(k^1,k^2)$, $k^1$ even, and $k^2$ odd or $k^1$ odd, and $k^2$ even.
\begin{multline*}
  \Nsymb^1(\boldsymbol{a})_{(-k^1,k^2)}=\sum_{k_1+k_2=k}{a^1_{(-k_1^1,k_1^2)}(-k_2^1)a^1_{(-k^1_2,k^2_2)}}+\sum_{k_1+k_2=k}{a^2_{(-k^1_1,k_1^2)}k^2_2a^1_{(-k^1_2,k^2_2)}}=\\
  +\sum_{k_1+k_2=k}{-(S\boldsymbol{a})^1_{(k_1^1,k_1^2)}(-k_2^1)(S\boldsymbol{a})^1_{(k^1_2,k^2_2)}}+\sum_{k_1+k_2=k}{(S\boldsymbol{a})^2_{(k^1_1,k_1^2)}k^2_2(S\boldsymbol{a})^1_{(k^1_2,k^2_2)}}=
\\\Nsymb^1(\symS\boldsymbol{a})_{(k^1,k^2)} = \Nsymb^1(\symS\boldsymbol{a})_{(-k^1,k^2)}.
\end{multline*}
If we consider indices $k_1,\,k_2$ such that $k=k_1+k_2$, it holds that both of the indices in one of the pairs $k^1_j,\,k^2_j$ are even (odd),
and in the second pair indices are of different parity (one even, and the other odd). This implies that in the last equality,
in the first term the symmetry generates single minus, as both of the modes come from the same component,
the second minus appears in front of the index $-k_2^1$. Whereas in the second term, as the modes come from different components,
the symmetry generates either two minuses ($k^1_1$ even, $k^2_1$ odd or vice-versa, and $k^1_2,\,k^2_2$ even (odd)), or
none minuses ($k^1_1$, $k^2_1$ even (odd), and $k^1_2$ even, $k^2_2$ odd or vice-versa). Finally, we obtain the claim.

\qed

\bigskip
We remark that there is another symmetry exhibited by the solutions of \eqref{eq:infDim}, which we denote by $\symS^{x\leftrightarrow y}$.
Existence of the second solution in Theorem~\ref{thm:main} follows from the bounds we establish in Section~\ref{sec:invOperatorBounds}
and the symmetry defined below.
\begin{definition}
  \label{defsymab}
\emph{Symmetry $\symS^{x\leftrightarrow y}$} by reflection by this symmetry we will obtain the existence of the second solution from Theorem~\ref{thm:main}.
Let $\symS^{x\leftrightarrow y}\colon H\times H\to H\times H$ be the following symmetry (denoted $x\leftrightarrow y$ in the prequel)
\begin{align*}
  \symS^{x\leftrightarrow y}\left(\mathbf{a}\right)^1_{(k^1,k^2)} &= \mathbf{a}^2_{(k^2,k^1)},\\
  \symS^{x\leftrightarrow y}\left(\mathbf{a}\right)^2_{(k^1,k^2)} &= \mathbf{a}^1_{(k^2,k^1)},\text{ for }k\in\mathbb{Z}.
\end{align*}
It is immediately verified that the solutions of the system \eqref{eq:toyInFourier} and all its Galerkin approximations are
invariant under this symmetry, i.e. $\boldsymbol{G}(\boldsymbol{a})=0\iff \boldsymbol{G}(\symS^{x\leftrightarrow y}\boldsymbol{a})=0$
as long as $\boldsymbol{F}$ is symmetric.
\end{definition}

\subsection{Structure of the linear operator}

Now, let us present the linear operator
\begin{equation}
\label{eq:lLambdaOper}
  L_\lambda(w)=\lambda\sin{y}\partial_xw+(-\Delta)^mw
\end{equation}
in Fouriers' coordinates introduced previously.
Here we argue how to reduce the problem of deriving dimension independent bounds for $L_\lambda^{-1}$
to the problem of bounding particular matrix norms. Recall that the operator $(-\Delta)^m$ in Fouriers' coordinates is diagonal
\begin{equation*}
  (-\Delta)^m={\rm diag\,}(\dots,l^{2m}+k^{2m},\dots).
\end{equation*}
In order to show the action of the $\lambda\sin{y}\partial_xw$ component, we introduce the following subspaces

\begin{definition}
  Let $l\in\{0,1,\dots,N\}$. We denote the following subspace of $H_{(N)}$
  \begin{equation*}
    H_{(N)}\supset H^l_{(N)} = \left\{w\colon w = \sum_{k=0}^{N}{w_k^{\pm}\sin\left(\pm lx + ky\right)} = \sum_{k=0}^{N}{-\frac12 a_{(\pm l,k)}\left(e^{i(\pm l, k)\cdot (x,y)} - e^{-i(\pm l, k)\cdot (x,y)}\right)},\ w_0^\pm = w_0\right\}.
  \end{equation*}
\end{definition}

It is easy to see that $H^l$ subspaces are invariant for the operator $L_\lambda(w)$ in the following sense
$  L_\lambda H^l_{(N)}\subset H^l_{(N+1)}$.

\begin{definition}
  Let us denote the projection of $L_\lambda$ onto $H^l_{(N)}$ by
\begin{equation*}
  L_\lambda^l:=P_{H^l_{(N)}}\circ L_\lambda\circ P_{H^l_{(N)}}\colon H^l_{(N)}\to H^l_{(N)}.
\end{equation*}
\end{definition}
Let
\[
  \symSpace_{(N)} = H^0_{(N)}\oplus H^1_{(N)}\oplus\cdots H^N_{(N)}.
\]
Define the projection of $L_\lambda$ onto the following space
\begin{definition}
  Let $R\subset H_{(N)}$ denote the following space
  \begin{equation}
    \label{subspaceR}
    R:=H^1\oplus\cdots\oplus H^N.
  \end{equation}
  The projection of $L_\lambda$ onto $R$ will be denoted by
  \begin{equation*}
    L_\lambda^R := P_{R}\left(L_\lambda \right)P_{R}\colon R\to R.
  \end{equation*}
\end{definition}

In order to present action of the operator $L_\lambda$ on a vector in $H^l_{(N)}$ we take
\[
  H^l_{(N)}\ni w=w_0\sin(lx)+\sum_{j=1}^{N}{w_j^\pm\sin(\pm lx+jy)}
\]

\begin{multline}
\label{eq:sinAction}
  \lambda\sin{y}w_{x}=l\lambda\sin{y}\left(w_0\cos{lx} + \sum_{j=1}^{N}{w_j^+\cos(lx+jy) - w_j^-\cos(-lx+jy)}\right)\\
  =\frac{l}{2}\lambda\left[ w_0\sin(lx+y) + w_0\sin(-lx+y) + \sum_{j = 1}^{N}{w_j^+\left[\sin(lx+(j+1)y) - \sin(lx+(j-1)y)\right]}\right.\\\nonumber
    \left. - \sum_{j = 1}^{N}{w_j^-\left[\sin(-lx+(j+1)y) - \sin(-lx+(j-1)y)\right]}\right],\nonumber
\end{multline}
Therefore
\begin{subequations}
\label{eq:triDiagSystem}
\begin{align}
  \left(L_{\lambda}^l w\right)_{l,0}&=-\frac{\lambda l}{2} w^+_{1}\ \ \, -\frac{\lambda l}{2} w^-_{1}\ \ +l^{2m}w_{0},\\
  \left(L_{\lambda}^l w\right)_{\pm l,j}&=\ \pm\frac{\lambda l}{2} w^\pm_{j-1}\, \mp\frac{\lambda l}{2} w^\pm_{j+1} + (l^{2m}+j^{2m})w^\pm_{j}\text{ for }j=1,\dots,N-1,\\
  \left(L_{\lambda}^l w\right)_{\pm N,j}&=\ \pm\frac{\lambda l}{2} w^\pm_{j-1} + (N^{2m}+j^{2m})w^\pm_{j},
\end{align}
\end{subequations}
where we used the convention $w_0^\pm = w_0$.

\providecommand{\symSpace}{H^\prime}

\label{sec:linearized}
We will study the structure of the linear operator $L_{\lambda}^l$ acting on the subspace $G^l_{(N)}$
  \begin{equation}
    \label{Gsubspace}      
    H^l_{(N)}\supset G^l_{(N)} = \left\{w\colon w\in\spann{\sin(lx),\sin(lx+ky)\text{ for }k=1,2,\dots,N}\right\}.
  \end{equation}  
  The subspace $G^l_{(N)}$ does not include in its span the part of the basis functions $\{\sin(-lx+ky),\ k=1,2,\dots,N\}$,
  which are present in the span of $H^l_{(N)}$. As we always work with vector solutions satisfying symmetry $S$ (Definition~\ref{defsyms}), for a given $(w_1,w_2)\in G^l_{(N)}\times G^l_{(N)}$ there is a unique $(v_1,v_2)\in H^l_{(N)}\times H^l_{(N)}$. In other words, the coefficients $(w_1^-, \dots, w_j^-, \dots, w_N^-)$ are determined
  by the corresponding coefficients with '+', i.e. $(w_1^+, \dots, w_j^+, \dots, w_N^+)$ through symmetry $S$ (Definition~\ref{defsyms}).

  In the sequel we will study the operator
\begin{equation*}
  L_\lambda^l:=P_{G^l_{(N)}}\circ L_\lambda\circ P_{G^l_{(N)}}\colon G^l_{(N)}\to G^l_{(N)},
\end{equation*}
which has the following tridiagonal form

\begin{align}
  \label{eq:Llambdal}
  L_\lambda^l=\left[\begin{array}{cccccc}
      l^{2m} &-\frac{l\lambda}{2}  &0          &         &\cdots   &0        \\
      \frac{l\lambda}{2}&l^{2m}+1  &-\frac{l\lambda}{2}  &0        &\cdots   &0        \\
      0      &\ddots    &\ddots     &\ddots   &\ddots   &0                 \\
      0      &\cdots    &0          &\frac{l\lambda}{2} &l^{2m} + (N-1)^{2m}&-\frac{l\lambda}{2}  \\
      0      &\cdots    &           &0        &\frac{l\lambda}{2}  &l^{2m} + N^{2m} \\
  \end{array}
  \right].
\end{align}
We will study the following full (projected) linear operator
\[
  P_{\oplus_{l=0}^{N} G^l_{(N)}}\circ L_\lambda\circ P_{\oplus_{l=0}^{N} G^l_{(N)}}\colon \oplus_{l=0}^{N} G^l_{(N)}\to \oplus_{l=0}^{N} G^l_{(N)}.
\]

In the sequel, we will use simply $L_\lambda$ to denote the full linear operator,
which has the following block diagonal form
\begin{equation}
  \label{eq:Llambda}
  L_\lambda=\left[\begin{array}{cccc}\operSS{0}&0&\cdots&0\\0&\operSS{1}&0&0\\\vdots&\ddots&\ddots&\vdots\\0&0&\operSS{N-1}&0\\0&\cdots&0&\operSS{N}\end{array}\right].
\end{equation}

\subsection{Bounds for matrices inverse to $L_\lambda^l$, $L_\lambda$ }
\label{sec:invOperatorBounds}
In this part we provide results on bounds of the particular norms of inverse tridiagonal matrices.
Some technical lemmas used to prove the presented bounds are provided in Section~\ref{sec:technicalLemmas}.
\begin{lemma}
 \label{lem:lInftyEstimate}
   Let $N>0$, $m>1$, $l=1,\dots,N$. The following uniform bound holds
\begin{equation*}
    \left|(\operS^{-1})_{j,k}\right|\leq \seriesAconst\left(\frac{l\lambda}{2}\right)^{-1},\quad\text{ for }k,j=1,\dots,N.
\end{equation*}
\end{lemma}

\begin{lemma}
  \label{lem:lOneEstimate}
  Let $N>0$, $m>1$, $l=1,\dots,N$. There exist $C(m)>0$ \emph{(independent of $\lambda$ and $N$)}, such that for
  $\lambda>2$ the following bounds hold
  \begin{equation*}
    \sum_{j=1,\dots,N}{\left|(\operS^{-1})_{ij}\right|}\leq C(m)\left(\frac{l\lambda}{2}\right)^{-1+1/2m},\qquad
    \sum_{j=1,\dots,N}{\left|(\operS^{-1})_{ji}\right|}\leq C(m)\left(\frac{l\lambda}{2}\right)^{-1+1/2m},
  \end{equation*}
  for all $i=1,\dots,N$.
\end{lemma}

In the next theorem we present the main result of this section, which is composed of bounds for the following norms 
$\|\operSS{R}^{-1}\|_{l^1\to l^1}$, $\|\operSS{R}^{-1}\|_{l^1\to l^\infty}$, $\|\operSS{R}^{-1}\|_{l^1\to l^1_1}$, see Definition~\ref{banach}. 
Where first two are standard norms, and the third (which we call the \emph{gradient norm}) is defined
\begin{definition}
Let $A\in\R^{N^2\times N^2}$ be a block-diagonal matrix
\begin{equation*}
  A=\left[\begin{array}{cccc}A_1&0&0&0\\\vdots&\ddots&\ddots&\vdots\\0&0&A_{N-1}&0\\0&\cdots&0&A_{N}\end{array}\right],
\end{equation*}
where $A_1,\dots,A_{N-1},A_N$ are $N$ dimensional square matrices.

We call the \emph{gradient norm} of $A$ the following matrix norm
  \begin{equation}
    \label{eq:gradientNorm}
    \|A\|_{l^1\to l^1_1}=\max_{j=1,\dots,N^2}{\left\|A(j)\right\|_{l^1_1}}=
    \max_{\substack{l=1,\dots,N\\j=1,\dots,N}}{\left\| A_{l}(j) \right\|_{l^1_1}}=
    \max_{\substack{l=1,\dots,N\\j=1,\dots,N}}{\sum_{k=1}^{N}{\left(l+k\right)\left|\left(A_{l}\right)_{k,j}\right|}},
  \end{equation}
\end{definition}
where $A(j)$ denotes the $j$-th column of $A$.

\begin{theorem}
\label{lem:invMatrixEstimates}
  Let $l=1,\dots,N$. Let $\operSS{l}$ be the matrix given by \eqref{eq:Llambdal}, $\operSS{R}$ be
the truncated matrix $P_RL_\lambda$ (projection of $L_\lambda$ onto the space $R$ \eqref{subspaceR}).

  The following estimates hold for the matrices $\operSS{l}^{-1}$ (diagonal submatrices of $\operSS{R}^{-1}$).
  \begin{align*}
     \|\operSS{l}^{-1}\|_{l^1\to l^1}&\leq C_1(m)\left(\frac{l\lambda}{2}\right)^{-1+1/2m},\\
     \|\operSS{l}^{-1}\|_{l^1\to l^\infty}&\leq\seriesAconst\left(\frac{l\lambda}{2}\right)^{-1},\\
     \|\operSS{l}^{-1}\|_{l^1\to l^1_1}&\leq C_2(m)\left(\frac{l\lambda}{2}\right)^{-1+1/m}.
   \end{align*}

  The following estimates hold for the matrix $\operSS{R}^{-1}$
   \begin{align*}
     \|\operSS{R}^{-1}\|_{l^1\to l^1}&\leq C_1(m) \left(\frac{\lambda}{2}\right)^{-1+1/2m},\\
     \|\operSS{R}^{-1}\|_{l^1\to l^\infty}&\leq\seriesAconst \left(\frac{\lambda}{2}\right)^{-1},\\
     \|\operSS{R}^{-1}\|_{l^1\to l^1_1}&\leq C_2(m) \left(\frac{\lambda}{2}\right)^{-1+1/m}.
   \end{align*}
\end{theorem}
We present a proof of this theorem in Section~\ref{sec:technicalLemmas}.

\providecommand{\A}{A}
\providecommand{\B}{B}
\providecommand{\Aa}{a}
\providecommand{\Bb}{b}
\providecommand{\AaBb}{\left(\Aa,\Bb\right)^T}

\section{Fixed point argument}
\label{sec:fixedPoint}

Having the estimate for the operator $(L^R_\lambda)^{-1}$ we are prepared to prove the  main result of the paper.
We assume that the considered solutions to \eqref{p2} are finite dimensional.
This assumption allows to use the results about the matrices norms presented in Section~\ref{sec:invOperatorBounds}.
Let us define two projections of the space $H$
\begin{align}\label{f4}
  P_R&\text{ projection onto }H^1\oplus\cdots\oplus H^N\text{ (the rotation like part) },\\ P_D&=I-P_R\text{ (the diagonal part) },
\end{align}
where $P_D$ is the projection onto the subspace free of $y$ dependence.
We proceed as follows. First, we construct an a-priori estimate for the solution of \eqref{eq:linearization}.  Let us display
basic features of $(A,B)\in\symSpace$ -- the solutions to \eqref{eq:linearization}, which follows directly from the bounds
presented in Theorem~\ref{lem:invMatrixEstimates} (where we absorb the $\frac12$ factor into the constant), namely
\begin{align*}
   \|\operSS{R}^{-1}\|_{l^1\to l^1}&\lesssim \lambda^{-1+1/2m},\\
   \|\operSS{R}^{-1}\|_{l^1\to l^\infty}&\lesssim \lambda^{-1},\\
   \|\operSS{R}^{-1}\|_{l^1\to l^1_1}&\lesssim \lambda^{-1+1/m}.
\end{align*}

Observe that due to the identities $P_DA,\,P_DB=0$ it is enough to use the bound for 
$(L^D_\lambda)^{-1}$, and we obtain
\begin{equation}\label{f1}
\begin{array}{lll}
  \|\B\|_{l^1}\lesssim \lambda^{1/2m}&\quad&\|\A\|_{l^1}\lesssim\lambda^{1/m},\\
  \|\B\|_{l^1_1}\lesssim \lambda^{1/m}&\quad&\|\A\|_{l^1_1}\lesssim\lambda^{3/2m},\\
  \|\B\|_{l^\infty} \lesssim 1 &\quad& \|\A\|_{l^\infty} \lesssim \lambda^{1/2m}.
\end{array}
\end{equation}
Consequently
\begin{align}\label{f2}
  \|(\A,\B)^T\cdot\nabla \A\|_{l^1}&\lesssim\|(\A,\B)^T\|_{l^1}\|\A\|_{l^1_1}\lesssim\lambda^{1/m}\lambda^{3/2m}\lesssim\lambda^{5/2m},\\
  \|(\A,\B)^T\cdot\nabla \B\|_{l^1}&\lesssim\|(\A,\B)^T\|_{l^1}\|\B\|_{l^1_1}\lesssim\lambda^{1/m}\lambda^{1/m}\,\lesssim\lambda^{2/m},\\
  \|(\A,\B)^T\cdot\nabla \A\|_{l^\infty}&\lesssim\|(\A,\B)^T\|_{l^\infty}\|\A\|_{l^1_1}\lesssim\lambda^{1/2m}\lambda^{3/2m}\lesssim\lambda^{2/m}.
\end{align}
In the estimations above, and generally in the estimates derived in this section we use often Young's inequality for products, i.e. 
\[
  \|f*g\|_1\leq\|f\|_1\|g\|_1\text{, and }\|f*g\|_\infty\leq\|f\|_\infty\|g\|_1.
\]

Now we split $V$ -- the solution to the system (\ref{p1}) in the following way
\begin{equation}\label{f3}
  V=\left(\begin{array}{c}\A\\\B\end{array}\right)+\left(\begin{array}{c}\Aa\\\Bb\end{array}\right),
\end{equation}
where $(A,B)^T$ is the solution of the linearized system \eqref{eq:linearization}.

In order to obtain the desired a priori estimate we are required to find a special property of function $b$. Namely, we prove that
\begin{equation}\label{Pb}
 P_D b =0,
\end{equation}
i.e. in $b$ there is no element depending only on $y$. To show $(\ref{Pb})$ we look at the rhs of (\ref{p2}) on the equation on $b$. 
We see that by (\ref{AB-con}) and (\ref{ab-con})
\begin{equation}
 (a+A)(b_x+B_x) + (b+B)(b_y+B_y)= (a+A)(b_x+B_x) + (b+B) (a_x+A_x) = \partial_x \left( (a+A)(b+B) \right).
\end{equation}
Hence 
\begin{equation}
 P_D \left( (a+A)(b_x+B_x) + (b+B)(b_y+B_y)\right) = 0, \mbox{ \ \ i.e. \ }  P_D b =0.
\end{equation}

{\bf Standing assumptions.} At the formal level of the a-priori estimate we assume that
 solutions to (\ref{p2}) fulfill
\begin{subequations}
  \label{eq:standingAssumptions}
\begin{align}
  &\left\|\AaBb\right\|_{l^1_1}\leq\lambda^{1-1/m},\label{eq:assumptionOne}\\
  &\|P_R\Aa\|_{l^1}\leq 1,\label{eq:assumptionTwo}\\
  &m>9/2.\label{eq:assumptionThree}
\end{align}
\end{subequations}
Recall (\ref{p2})
\begin{equation}\label{f6}
 \begin{array}{l}
  \lambda \sin y \d_x a + (-\Delta)^m a = - \lambda \cos y \,  b - \left(
  \begin{array}{c}
   a + A \\
   b + B
  \end{array}
  \right) \cdot \nabla (a+A), \\
  \lambda \sin y \d_x b + (-\Delta)^m b =  - \left(
  \begin{array}{c}
   a + A \\
   b + B
  \end{array}
  \right) \cdot \nabla (b+B).
 \end{array}
\end{equation}

In order to find the bound we apply the estimates for $L^R_\lambda$ formally, assuming that the solutions are finite dimensional.
Treating the right hand side of \eqref{f6} we have the following bounds

\paragraph{Bound for $\|\Bb\|_{l^1_1}=\|P_R\Bb\|_{l^1_1}$.}
\begin{align}
  \|\Bb\|_{l^1_1}&\lesssim\|{L^R_\lambda}^{-1}\|_{l^1\to l_1^1}\left(\|(\A,\B)^T\cdot\nabla\B\|_{l^1}+\|\Aa\Bb_x\|_{l^1}+\|\Bb\Bb_y\|_{l^1}
  +\|\Aa\B_x\|_{l^1}+\|\Bb\B_y\|_{l^1}+\|(\A,\B)^T\cdot\nabla\Bb\|_{l^1}\right)\nonumber\\
&\lesssim\lambda^{-1+1/m}\left(\lambda^{2/m}+\|\Aa\|_{l^1}\|\Bb\|_{l^1_1}+\|\Bb\|_{l^1}\|\Bb\|_{l^1_1}+\|\Aa\|_{l^1}\lambda^{1/m}+\|\Bb\|_{l^1}\lambda^{1/m}+\lambda^{1/m}\|\Bb\|_{l^1_1}\right)\nonumber\\
&\lesssim \lambda^{-1+3/m}+\|\Aa\|_{l^1}\lambda^{-1+2/m},\label{eq:boundForB}
\end{align}
where the last inequality is obtained after cleaning the absorbed terms, which is due to the assumptions \eqref{eq:assumptionOne}, and
\eqref{eq:assumptionThree}.
We will also need the following estimate for $\|\Bb\|_{l^1}$, derived analogously as above
\begin{equation}
  \label{eq:blOneEstimate}
  \|\Bb\|_{l^1}\lesssim \lambda^{-1+5/2m}+\|\Aa\|_{l^1}\lambda^{-1+3/2m}.
\end{equation}

Let us define
\begin{equation}\label{f8}
  \|P_D\Aa\|_{l^\infty_{2m}}=\sup_{\substack{k\in\mathbb{Z}\\|k|\leq N}}{k^{2m}|a_{(0,k)}|}.
\end{equation}
\paragraph{Bound for $\|P_D\Aa\|_{l^\infty_{2m}}$.}
In this case the operator $P_DL_\lambda$ is diagonal, therefore we bound the particular norm $\|P_D\Aa\|_{l^\infty_{2m}}$, it is trivially
bounded by the $l^\infty$ norm of the right hand side. Moreover, observe that $l^\infty_{2m}$ norm bounds $l^1_1$, i.e. we have
$\|P_D\Aa\|_{l^1_1}\lesssim\|P_D\Aa\|_{l^\infty_{2m}}$ for $m>3/2$, remembering that the dimension is two.
\begin{align}
  \|P_D\Aa\|_{l^\infty_{2m}}&\lesssim\|(\A,\B)^T\cdot\nabla A\|_{l^\infty}+\|P_D(\Aa\Aa_x) \|_{l^\infty}+\|P_D(\Bb\Aa_y)\|_{l^\infty}+
\|P_D(\Aa\A_x)\|_{l^\infty}+\|P_D(\Bb\A_y) \|_{l^\infty}\nonumber\\
& +\|P_D((\A,\B)^T\cdot\nabla\Aa)\|_{l^\infty}\nonumber\\
&\lesssim\lambda^{2/m}+\|P_R\Aa\|_{l^1}\|P_R\Aa\|_{l^1_1}+\|\Bb\|_{l^1}\|P_R\Aa\|_{l^1_1}
+\|P_R\Aa\|_{l^1}\lambda^{3/2m}+\|\Bb\|_{l^1}\lambda^{3/2m}+\lambda^{1/m}\|P_R\Aa\|_{l^1_1}.\label{eq:bbound}
\end{align}
We removed all terms, which do not generate $P_D$, i.e. any product of terms, one of them being in $P_D$, and the other one in $P_R$.
When the bound \eqref{eq:bbound} is used (potentially the worst term $\|P_D\Aa P_D\Aa_x\|$ is not present as $P_D\Aa_x=0$) we get
\begin{align*}
  \|P_D\Aa\|_{l^\infty_{2m}}\lesssim&\lambda^{2/m}+\|P_R\Aa\|^2_{l^1_1}+\lambda^{-1+3/m}\|P_R\Aa\|_{l^1_1}+\|\Aa\|_{l^1}\|P_R\Aa\|_{l^1_1}\lambda^{-1+2/m}\\
&+\|P_R\Aa\|_{l^1}\lambda^{3/2m}+\lambda^{-1+9/2m}+\|\Aa\|_{l^1}\lambda^{-1+7/2m}+\lambda^{1/m}\|P_R\Aa\|_{l^1_1}\\
\lesssim&\lambda^{2/m}.
\end{align*}

To get last inequality we used the assumption \eqref{eq:assumptionTwo} and \eqref{eq:assumptionThree}, the term $\lambda^{2/m}$
is clearly of the highest order from the terms that are left. Here we use that $m > 9/2$.

\paragraph{Bound for $\|P_R\Aa\|_{l^1_1}$.}

Observe that we have $\|P_D\Aa\|_{l^1_1}\lesssim\|P_D\Aa\|_{l^\infty_{2m}}\lesssim\lambda^{2/m}$
\begin{align*}
  \|P_R\Aa\|_{l^1_1}\lesssim&\lambda^{-1+1/m}\left(\lambda\|\Bb\|_{l^1}+\|(\A,\B)^T\cdot\nabla\A\|_{l^1}+\|\Aa\|_{l^1}\|\Aa\|_{l^1_1}+\|\Bb\|_{l^1}\|\Aa\|_{l^1_1}
+\|\Aa\|_{l^1}\|\A\|_{l^1_1} \right. \\
& \left. +\|\Bb\|_{l^1}\|\A\|_{l^1_1}+\|(\A,\B)^T\|_{l^1}\|\Aa\|_{l^1_1}\right)\\
\lesssim & \lambda^{-1+1/m}\left(\lambda\|\Bb\|_{l^1}+\lambda^{5/2m}+\|P_R\Aa\|^2_{l^1_1}+\|P_R\Aa\|_{l^1_1}\lambda^{2/m}+
\|\Bb\|_{l^1}\|P_R\Aa\|_{l^1_1}+\right.\\
&\left.\lambda^{2/m}\lambda^{3/2m} + \|P_Ra\|_{l^1_1}\lambda^{3/2m} + \|\Bb\|_{l^1}\lambda^{3/2m} +
\lambda^{1/m} \lambda^{2/m}\right)\\
\lesssim&\lambda^{-1+1/m}\left(\lambda\|\Bb\|_{l^1}+\|P_R\Aa\|^2_{l^1_1}+\|P_R\Aa\|_{l^1_1}\lambda^{2/m}+\|\Bb\|_{l^1}\lambda^{3/2m}+
\|\Bb\|_{l^1}\|P_R\Aa\|_{l^1_1}+\lambda^{7/2m}\right).
\end{align*}
Observe that after the second inequality the term $\|b\|_{l^1}\|a\|_{l^1_1}$ is not present as $P_Da_x=0$,
clearly the highest order term is $\|P_Da\|_{l^1}\|A\|_{l^1_1}=\lambda^{2/m}\lambda^{3/2m}=\lambda^{7/2m}$.

Now we use the bound \eqref{eq:bbound}, and remove some of the terms that were absorbed by using the assumptions
\eqref{eq:assumptionOne}, \eqref{eq:assumptionTwo}, and \eqref{eq:assumptionThree}, observe in the inequality above 
the bad looking term $\lambda\|\Bb\|_{l^1}$, we estimate it using \eqref{eq:blOneEstimate}
\begin{equation*}
  \lambda\|\Bb\|_{l^1}\lesssim \lambda^{5/2m}+\lambda^{3/2m}\left(\|P_D\Aa\|_{l^1}+\|P_R\Aa\|_{l^1}\right),
\end{equation*}
\begin{multline*}
  \|P_R\Aa\|_{l^1_1}\lesssim\lambda^{-1+1/m}\left(\lambda^{5/2m}+\lambda^{2/m}\lambda^{3/2m}+\|P_R\Aa\|^2_{l^1_1}+\|P_R\Aa\|_{l^1_1}\lambda^{2/m}+\right.\\
\left.+(\lambda^{-1+5/2m}+\|\Aa\|_{l^1}\lambda^{-1+3/2m})(\lambda^{3/2m}+\|P_R\Aa\|_{l^1_1})+\lambda^{7/2m}\right).
\end{multline*}
After using the assumption \eqref{eq:assumptionTwo} all terms with $\|P_R\Aa\|_{l^1_1}$ are being absorbed,
and clearly the highest order term in the parenthesis is $\lambda^{7/2m}$, so finally we end up with
\begin{align}\label{f10}
  \|P_R\Aa\|_{l^1_1}&\lesssim\lambda^{-1+9/2m}.
\end{align}
Observe that $P_R\Aa$ is mapped into itself by the operator $L_\lambda^{-1}$, due to the assumption \eqref{eq:assumptionThree},
namely $m > 9/2$.

Going back to \eqref{eq:boundForB} we get that
\begin{equation}\label{f11}
  \|\Bb\|_{l^1_1}\leq \lambda^{-1+4/m}.
\end{equation}

Summing up the considerations from this part we obtain the following result

\smallskip

\begin{lemma}
\label{lem:apriori}
Let $a,b$ be a small solution to problem \eqref{f6}, then it obeys the following \emph{dimension independent a-priori estimate}
\begin{equation}\label{f12}
 \|b\|_{l^1_1} \leq C \lambda^{-1+4/m}, \qquad \| P_R a\|_{l^1_1}\leq C\lambda^{-1+9/2m}, \qquad \|P_D a\|_{l^1_1}\leq C \lambda^{2/m}.
\end{equation}

\end{lemma}

\section{Proof of main theorem}
\label{sec:proof}
Using the so far presented results, we may now proceed to proving our main result -- Theorem~\ref{thm:main}.
Here we want to construct the solutions, using the system \eqref{f6} and the a-priori estimates (Lemma~\ref{lem:apriori}).

We start with the construction of the sequence of solution's approximations.
We define the solution $(a_{n+1},b_{n+1})$ as the solution to the following problem

\begin{equation}\label{f13}
 \begin{array}{l}
  \lambda \sin y \d_x a_{n+1} + (-\Delta)^m a_{n+1} = - \lambda \cos y \,  b_{n+1} - \left(
  \begin{array}{c}
   a_n + A_n \\
   b_n + B_n
  \end{array}
  \right) \cdot \nabla (a_n+A_n), \\
  \lambda \sin y \d_x b_{n+1} + (-\Delta)^m b_{n+1} =  - \left(
  \begin{array}{c}
   a_n + A_n \\
   b_n + B_n
  \end{array}
  \right) \cdot \nabla (b_n+B_n).
 \end{array}
\end{equation}
We take $(a_0,b_0)=(0,0)$ and define $A_n,B_n$ as the projection onto the spaces $H_{(N+2)}$, where $N$ determines the number of active modes. If $(a_n,b_n) \in H_{(N)}$, then $(a_{n+1},b_{n+1}) \in H_{(2N+2)}$.
 
Note, in addition, that (\ref{f13}) guarantees the constraint (\ref{ab-con}). It is clear that from $\partial_x a_n = \partial_y b_n$ it follows
\begin{equation}\label{f13a}
 \partial_x \left[ \left(
  \begin{array}{c}
   a_n + A_n \\
   b_n + B_n
  \end{array}
  \right) \cdot \nabla (a_n+A_n) \right] - 
   \partial_y \left[\left(
  \begin{array}{c}
   a_n + A_n \\
   b_n + B_n
  \end{array}
  \right) \cdot \nabla (b_n+B_n)\right] \equiv 0.
\end{equation}
And this implies $\partial_x a_{n+1} = \partial_y b_{n+1}$. Thus constraint (\ref{ab-con}) is guaranteed, refer (\ref{p2b}).

Repeating the estimates for the system \eqref{f6} we find that if
\begin{equation}\label{f14}
 \|b_n\|_{l^1_1} \leq C \lambda^{-1+4/m}, \qquad \| P_R a_n\|_{l^1_1}\leq C\lambda^{-1+9/2m}, \qquad \|P_D a_n\|_{l^1_1}\leq C \lambda^{2/m}.
\end{equation}
then
\begin{equation}\label{f15}
 \|b_{n+1}\|_{l^1_1} \leq C \lambda^{-1+4/m}, \qquad \| P_R a_{n+1}\|_{l^1_1}\leq C\lambda^{-1+9/2m}, \qquad \|P_D a_{n+1}\|_{l^1_1}\leq C \lambda^{2/m}.
\end{equation}
with the same constants $C$, provided $\lambda$ sufficiently large.

We shall underline that for a fixed $n$ we are allowed to apply results for the finite dimensional approximation of $L_\lambda$.
We emphasize that all constants in Theorem~\ref{lem:invMatrixEstimates} are  independent on $N$.

We want to prove that $\{a_n,b_n\}$ is a Cauchy sequence. We consider the following system
\begin{equation}\label{f16}
 \begin{array}{l}
  \lambda \sin y \d_x (a_{n+1}-a_n) + (-\Delta)^m (a_{n+1}-a_n) = - \lambda \cos y \,  (b_{n+1}-b_n)
\\
\qquad \qquad - \left(
  \begin{array}{c}
   a_n + A_n \\
   b_n + B_n
  \end{array}
  \right) \cdot \nabla (a_n+A_n)
+\left(
  \begin{array}{c}
   a_{n-1} + A_{n-1} \\
   b_{n-1} + B_{n-1}
  \end{array}
  \right) \cdot \nabla (a_{n-1}+A_{n-1}), \\[20pt]
  \lambda \sin y \d_x (b_{n+1}-b_n) + (-\Delta)^m (b_{n+1} -b_n)=
\\
\qquad \qquad  - \left(
  \begin{array}{c}
   a_n + A_n \\
   b_n + B_n
  \end{array}
  \right) \cdot \nabla (b_n+B_n)
+ \left(
  \begin{array}{c}
   a_{n-1} + A_{n-1} \\
   b_{n-1} + B_{n-1}
  \end{array}
  \right) \cdot \nabla (b_{n-1}+B_{n-1}).
 \end{array}
\end{equation}
Taking a large $n$ we want to prove that
\begin{multline}\label{f17}
 \|P_R (a_{n+1}-a_n),b_{n+1}-b_n\|_{l^1_1}+\lambda^{-2/m}\|P_D(a_{n+1}-a_n)\|_{l^1_1} \leq\\\frac12
\left(\| P_R (a_{n}-a_{n-1}),b_{n}-b_{n-1}\|_{l^1_1} +\lambda^{-2/m}\|P_D(a_n-a_{n-1})\|_{l^1_1}\right)+ \epsilon_n,
\end{multline}
where $\epsilon_n \to 0$ as $n \to \infty$, the quantity $\epsilon_n$ is related by norms of terms like $(A_n-A_{n-1})$ and $(B_n-B_{n-1})$.

In order to justify (\ref{f17}) we point out few estimates which provides the inequality. Here we use the same tools as in the proof of Lemma \ref{lem:apriori}. 
Hence we estimate the right hand side of (\ref{f16}).
 We have to estimate the following terms.

\smallskip 

\underline{For $\|b_{n+1}-b_{n}\|_{l^1_1}$} we have
\begin{equation}
\|
\left(
\begin{array}{c}
 a_n + A_n \\
b_n + B_n
\end{array}
\right)
\cdot \nabla (B_n-B_{n-1})
\|_{l^1} \leq
\epsilon_n.
\end{equation}
For $n$ sufficiently large it it clear that $\|B_n-B_{n-1}\|_{l^1_1} \to 0$ as $n \to \infty$.
Next,
\begin{equation}
 \|
\left(
\begin{array}{c}
 a_n + A_n \\
b_n + B_n
\end{array}
\right)
\cdot \nabla  (b_n-b_{n-1})
\|_{l^1} \lesssim
\lambda^{2/m} \|b_n-b_{n-1}\|_{l^1_1}
\end{equation}
and
\begin{multline}
 \|
\left(
\begin{array}{c}
 a_n -a_{n-1}+ A_n-A_{n-1} \\
b_n -b_{n-1}+ B_n -B_{n-1}
\end{array}
\right)
\cdot \nabla (b_{n-1}+B_{n-1})
\|_{l^1} \\
\lesssim
\epsilon_n + (\lambda^{-1+4/m} +\lambda^{1/m})\|(a_n -a_{n-1}), b_n-b_{n-1}\|_{l^1_1}\\
\lesssim \epsilon_n + (\lambda^{-1+4/m} +\lambda^{1/m})\|P_R(a_n -a_{n-1}), b_n-b_{n-1}\|_{l^1_1}\\
+(\lambda^{-1+4/m} +\lambda^{1/m})\lambda^{2/m}\lambda^{-2/m}\|P_D(a_n -a_{n-1})\|_{l^1_1}.
\end{multline}

Hence
\begin{multline}\label{x1}
 \|b_{n+1}-b_n\|_{l^1_1} \lesssim \lambda^{-1 + 3/m} \|P_1(a_n -a_{n-1}), b_n-b_{n-1}\|_{l^1_1} \\
+ (\lambda^{-2+7/m}+\lambda^{-1+4/m})
\lambda^{-2/m}\|P_D(a_n-a_{n-1})\|_{l^1_1} + \epsilon_n.
\end{multline}

\smallskip 

\underline{For $\|P_D(a_{n+1}-a_n)\|_{l^1_1}$} we have
\begin{equation}
 \|
P_D\left(
\left(
\begin{array}{c}
 a_n + A_n\\
b_n+B_n
\end{array}
\right)
\cdot \nabla (a_n-a_{n-1})
\right)
\|_{l^1}
\lesssim \lambda^{1/m} \|P_R(a_n-a_{n-1})\|_{l^1_1},
\end{equation}
\begin{equation}
 \| P_D\left(\left( 
\begin{array}{c}
 a_n-a_{n-1} \\
b_n - b_{n-1}
\end{array}
\right)
\nabla A_n\right)\|_{l_1} \lesssim \lambda^{3/2m} \|P_R(a_n -a_{n-1}), b_n -b_{n-1}\|_{l^1_1}.
\end{equation}
The remaining terms here are simpler.
So
\begin{equation}\label{x2}
 \|P_D(a_{n+1}-a_n)\|_{l^1_1} \lesssim \lambda^{3/2m} \|P_R(a_n-a_{n-1}),b_n-b_{n-1}\|_{l^1_1} + \mbox{ better terms}.
\end{equation}

\smallskip 

\underline{ For $\|P_R(a_{n+1}-a_n)\|_{l^1_1}$} we have
\begin{multline}
 \|
P_R\left(
\left(
\begin{array}{c}
 a_n + A_n\\
b_n+B_n
\end{array}
\right)
\cdot \nabla (a_n-a_{n-1})
\right)
\|_{l^1}
\\
\lesssim \lambda^{2/m} \|P_R(a_n-a_{n-1})\|_{l^1_1} + (\lambda^{-1+4/m}+\lambda^{1/m})\lambda^{2/m}\lambda^{-2/m}\|P_D(a_n-a_{n-1})\|_{l^1_1}
\end{multline}
and
\begin{equation}
 \|P_R\left(\left(
\begin{array}{c}
 a_n-a_{n-1}\\
b_n-b_{n-1}
\end{array}
\right)\cdot
\nabla A_n\right)\|_{l^1} \lesssim \lambda^{3/2m}\|P_R(a_n-a_{n-1}),b_n-b_{n-1}\|_{l^1}+\lambda^{7/2m}\lambda^{-2/m}\|P_D(a_n-a_{n-1})\|_{l^1}.
\end{equation}

The last term is $\|\lambda \cos y (b_{n+1}-b_n)\|_{l^1}$, and using the estimates  for $\|L^B_\lambda\|_{l^1\to l^1}$ 
from Theorem~\ref{lem:invMatrixEstimates} we find
\begin{multline}
 \|\lambda \cos y (b_{n+1} - b_n)\|_{l^1} \lesssim \lambda^{1/2m} \|RHS (\ref{f16})_2\|_{l^1} \\
\lesssim
\lambda^{1/2m} \left(\lambda^{2/m} \|P_R(a_n -a_{n-1}), b_n-b_{n-1}\|_{l^1_1} \right.\\\left.
+ \lambda^{3/m}\lambda^{-2/m}\|P_D(a_n-a_{n-1})\|_{l^1_1} + \epsilon_n\right).
\end{multline}
Hence
\begin{align}\label{x3}
 \|P_R(a_{n+1} -a_n)\|_{l^1_1}&\\\lesssim&\lambda^{-1+1/m}\left(\lambda^{5/2m} \|P_R(a_n -a_{n-1}), b_n-b_{n-1}\|_{l^1_1}+
 \lambda^{7/2m}\lambda^{-2/m}\|P_D(a_n-a_{n-1})\|_{l^1_1}+\epsilon_n\right)\\
\lesssim&\lambda^{-1+7/2m}\|P_R(a_n-a_{n-1}),b_n-b_{n-1}\|_{l^1_1}+ \lambda^{-1+9/2m}\lambda^{-2/m}\|P_D(a_n-a_{n-1})\|_{l^1_1} + \epsilon_n.
\end{align}

\smallskip 

 Summing up (\ref{x1}),(\ref{x2}) and (\ref{x3}) we conclude
\begin{multline}
  \|P_R (a_{n+1}-a_n),b_{n+1}-b_n\|_{l^1_1}+\lambda^{-2/m}\|P_D(a_n-a_{n-1})\|_{l^1_1}
\lesssim
\\
\lambda^{-1+7/2m}\| P_R (a_{n}-a_{n-1}),b_{n}-b_{n-1}\|_{l^1_1} + \lambda^{-1+9/2m}\lambda^{-2/m}\|P_D(a_n-a_{n-1})\|_{l^1_1} + \epsilon_n
\end{multline}
For $m>9/2$ and large $\lambda$ we got (\ref{f17}).

The condition (\ref{f17}) implies that the sequence $\{P_Ra_n,b_n,P_Da_n\}$ has a limit in the space $l^1_1$. It means that there exists a solution
to problem (\ref{f6}) obeying estimates from Lemma \ref{lem:apriori}. In other words we have constructed the solution (\ref{i5}).
 We shall underline that the limit in $l^1_1$ implies that the derivative is uniformly bounded, 
thus the nonlinear term is described pointwisely. A bootstrap method implies that the solutions constructed in the above way are indeed smooth. 

\smallskip 

Existence of the solution (\ref{i6}) follows from the symmetry $x \leftrightarrow y$ from Definition~\ref{defsymab} that we recall here for completeness
\begin{align*}
  \symS^{x\leftrightarrow y}\left((a,b)\right)^1_{(k^1,k^2)} &= b_{(k^2,k^1)},\\
  \symS^{x\leftrightarrow y}\left((a,b)\right)^2_{(k^1,k^2)} &= a_{(k^2,k^1)},
\end{align*}
where it should be understood on the level of the Fourier modes of $(a,b)$.
Largeness of $\lambda$ implies there are two different solutions. Theorem \ref{thm:main} is proved.

\section{Analysis of large matrices and proof of Theorem \ref{lem:invMatrixEstimates}}
\label{sec:technicalLemmas}
\paragraph{Notation}
Let $N>0$ be an even number, $m>1$, $l>0$, $\lambda\in\mathbb{R}$.\\
Let us denote
\begin{equation*}
  \mathbb{R}^{2\times 2}\ni \matrixS^l(a,b)=\left[\begin{array}{cc}a&-\frac{l\lambda}{2}\\\frac{l\lambda}{2}&b\end{array}\right].
\end{equation*}
We denote a tridiagonal matrix with elements $\{a_j\}_{j=1}^N$ on the diagonal, $-l\lambda$ over diagonal, and $l\lambda$ under diagonal by
\begin{equation*}
  \mathbb{R}^{N\times N}\ni \matrixS^l(a_1,\dots, a_N)=\left[
  \begin{array}{cccccc}
    a_1&-\frac{l\lambda}{2}&0&0&0&\dots\\
    \frac{l\lambda}{2}&a_2&-\frac{l\lambda}{2}&0&0&\dots\\
    \ &\ &\ddots&\ &\ddots&\ \\
    \ &\ &\ &0&\frac{l\lambda}{2}&a_N\\
  \end{array}\right].
\end{equation*}
Let the increasing sequence $\{\seriesD[l]{j}\}_{j=1}^N$ be given by
\begin{equation*}
    \seriesD[l]{1}=l^{2m},\ \seriesD[l]{2}=l^{2m}+1,\cdots,\ \seriesD[l]{N}=l^{2m}+(N-1)^{2m}\quad\text{ for }l>0.
\end{equation*}
We denote the tridiagonal matrix with the increasing sequence $\{\seriesD[l]{j}\}_{j=1}^N$ on the diagonal by
\begin{equation*}
  \operS:=\matrixS^l(\seriesD[l]{1},\dots, \seriesD[l]{N})\in\mathbb{R}^{N\times N}.
\end{equation*}

\begin{figure}[htbp]
  \centering
        \begin{subfigure}[b]{0.5\textwidth}
          \includegraphics[width=\textwidth]{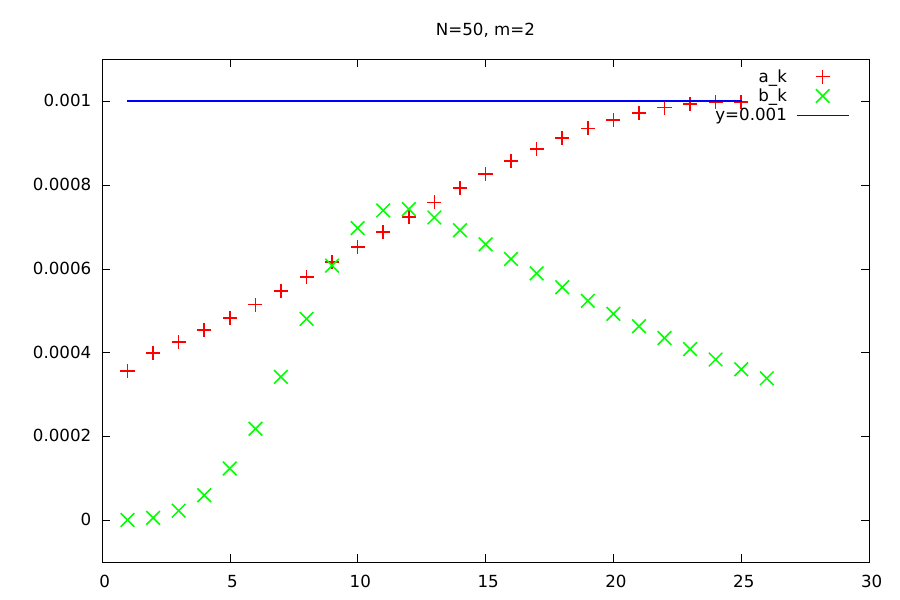}
        \end{subfigure}
        \begin{subfigure}[b]{0.5\textwidth}
          \includegraphics[width=\textwidth]{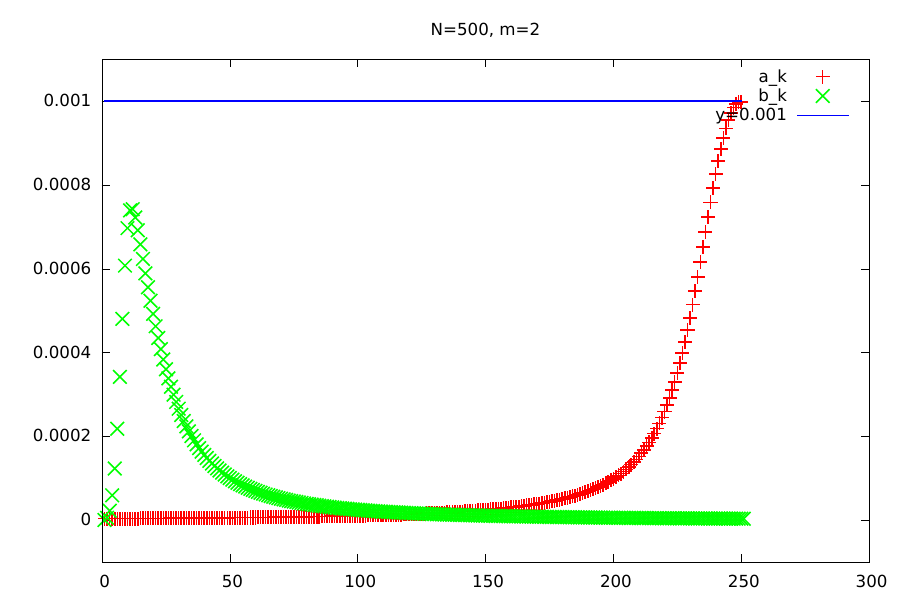}
        \end{subfigure}
        \begin{subfigure}[b]{0.5\textwidth}
          \includegraphics[width=\textwidth]{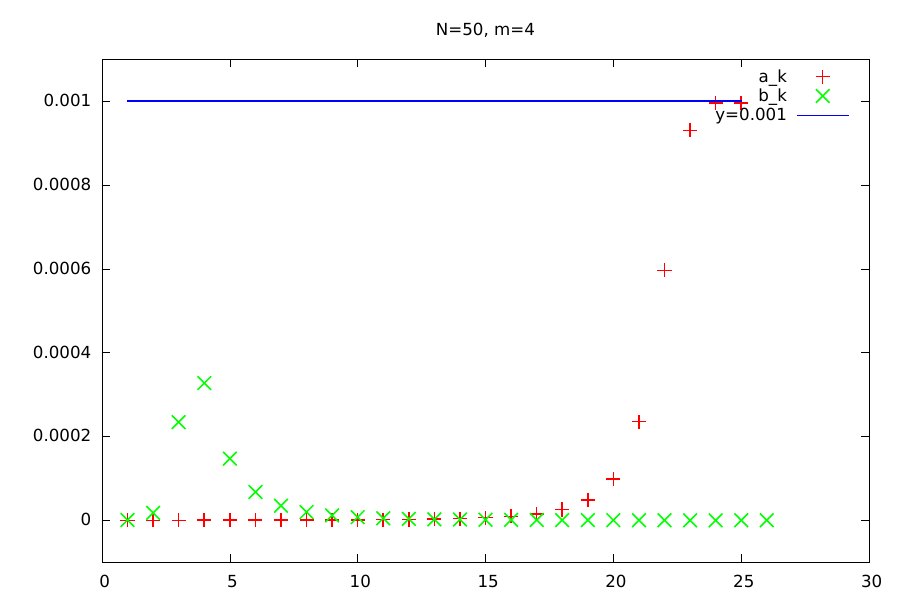}
        \end{subfigure}
        \begin{subfigure}[b]{0.5\textwidth}
          \includegraphics[width=\textwidth]{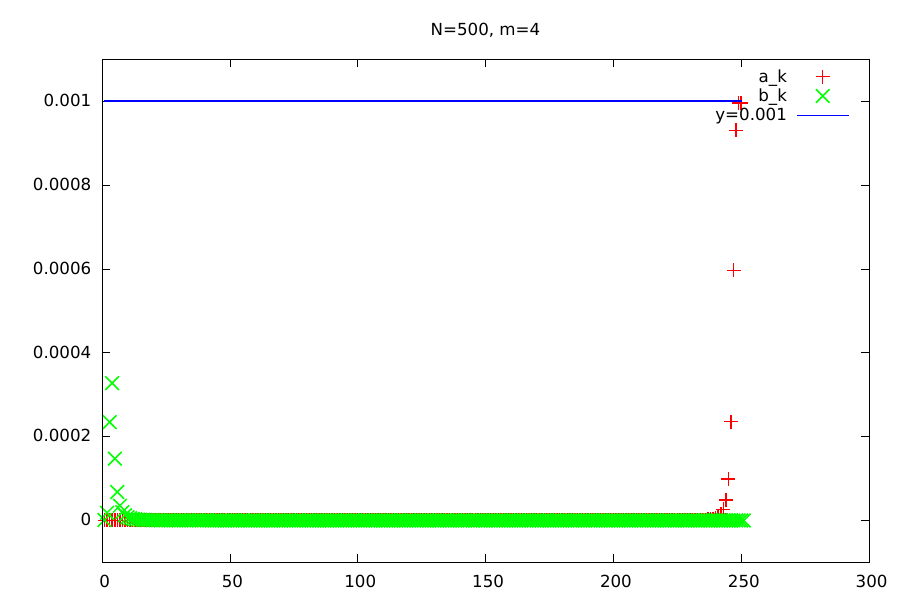}
        \end{subfigure}
  \caption{Graph showing numerically calculated recursive series $\{\seriesA[1]{j}\}$, $\{\seriesB[1]{j}\}$ for $\lambda=1000$, and for $m=2$ (top figures), $m=4$ (bottom figures). Apparent upper bound for the series is also shown $(1/\lambda)$.}
  \label{fig:series}
\end{figure}

\begin{definition}
  In the sequel we will use the following notation to denote the off diagonal term of the tridiagonal matrices $T^l$
  \[
    \hat{\lambda} = \hat{\lambda}(l) =  \frac{l\lambda}{2}.
  \] 
\end{definition}

\begin{lemma}
  \label{lem:sequences}
  Let $l>1$. Let the sequences $\{\seriesA[l]{j}\}$, $\{\seriesB[l]{j}\}$ be given by the following recursive formulas
  \begin{align*}
    &\seriesA[l]{0}=0,&\quad&\seriesB[l]{0}=0,\\
    &\seriesA[l]{1}=\frac{\seriesD[l]{N}}{\seriesD[l]{N-1}\seriesD[l]{N} + \hat{\lambda}^2},&\quad&\seriesB[l]{1}=\frac{\seriesD[l]{1}}{\seriesD[l]{1}\seriesD[l]{2}+\hat{\lambda}^2},\\
    &\seriesA[l]{2}=\frac{\seriesD[l]{N-2}+\seriesA[l]{1}\hat{\lambda}^2}{\seriesD[l]{N-3}\seriesD[l]{N-2}+\seriesA[l]{1}\seriesD[l]{N-3}\hat{\lambda}^2+\hat{\lambda}^2},&\quad&\seriesB[l]{2}=\frac{\seriesD[l]{3}+\seriesB[l]{1}\hat{\lambda}^2}{\seriesD[l]{3}\seriesD[l]{4}+\seriesB[l]{1}\seriesD[l]{4}\hat{\lambda}^2+\hat{\lambda}^2},\\
    &\seriesA[l]{j}=\frac{\seriesD[l]{N-2j+2}+\seriesA[l]{j-1}\hat{\lambda}^2}{\seriesD[l]{N-2j+1}\seriesD[l]{N-2j+2}+\seriesA[l]{j-1}\seriesD[l]{N-2j+1}\hat{\lambda}^2+\hat{\lambda}^2},&\quad&\seriesB[l]{j}=\frac{\seriesD[l]{2j-1}+\seriesB[l]{j-1}\hat{\lambda}^2}{\seriesD[l]{2j}\seriesD[l]{2j-1}+\seriesB[l]{j-1}\seriesD[l]{2j}\hat{\lambda}^2+\hat{\lambda}^2},
  \end{align*}
  for $j\geq 1$. Then the following bounds hold for all $j\geq 0$
  \begin{equation*}
    0\leq \seriesA[l]{j}\leq \seriesAconst / \hat{\lambda} ,\quad 0\leq\seriesB[l]{j}\leq 1 / \hat{\lambda}.
  \end{equation*}

\end{lemma}
\paragraph{Proof}
The $\seriesA[l]{j}, \seriesB[l]{j}\geq 0$ part of the bound is trivial.

Now we prove $\seriesAconst\hat{\lambda}^{-1}\geq \seriesA[l]{j}$. We proceed by induction, first we prove that
$\seriesAconst\hat{\lambda}^{-1}\geq \seriesA[l]{1}$ holds.
Observe that for all $N>1$, and $l\geq 1$ we have
\begin{multline}
  \label{eq:abestimate}
  \frac{l^{2m}+(N-1)^{2m}}{(l^{2m}+(N-1)^{2m})(l^{2m}+(N-2)^{2m})+\hat{\lambda}^2}<
  \frac{l^{2m}+2^{2m-1}\left[(N-2)^{2m}+1\right]}{(l^{2m}+(N-2)^{2m})^2+\hat{\lambda}^2}<\\
  \frac{2^{2m}\left[l^{2m}+(N-2)^{2m}\right]}{(l^{2m}+(N-2)^{2m})^2+\hat{\lambda}^2}<
  \frac{2^{2m}a}{a^{2}+\hat{\lambda}^2} = \frac{2^{2m}}{a+\frac{\hat{\lambda}^2}{a}}<\frac{2^{2m}}{\hat{\lambda}},
\end{multline}
where $a=l^{2m}+(N-2)^{2m}$, we used the estimate due to convexity $(N-1)^{2m}<2^{2m}\left(\frac{(N-2)+1}{2}\right)^{2m}<
2^{2m-1}\left({(N-2)^{2m}+1}\right)$, the last inequality follows from $a^2+\hat{\lambda}^2>a\hat{\lambda}$.

Assuming $\seriesAconst\hat{\lambda}^{-1}\geq a_{j-1}$ we verify that $\seriesAconst\hat{\lambda}^{-1}\geq a_{j}$ holds.

First, observe that $f(a):=\frac{\seriesD{N-2j+2}+a\hat{\lambda}^2}{\seriesD{N-2j+2}\seriesD{N-2j+1}+a \hat{\lambda}^2\seriesD{N-2j+1}+\hat{\lambda}^2}$ is a strictly increasing function
for all $a\geq 0$ (denominator is positive), as
\begin{equation}
  \label{eq:strictlyIncreasing}
  f^\prime(a)=\frac{\hat{\lambda}^4}{\left(\seriesD{N-2j+2}\seriesD{N-2j+1}+a \hat{\lambda}^2\seriesD{N-2j+1}+\hat{\lambda}^2\right)^2}\geq 0,
\end{equation}
so we have
\begin{equation}
  \label{eq:seriesAbound}
  \frac{\seriesD{k+2}+a_{j-1}\hat{\lambda}^2}{\seriesD{k+2}\seriesD{k+1}+a_{j-1}\seriesD{k+1}\hat{\lambda}^2+\hat{\lambda}^2}\leq
  \frac{l^{2m}+(k+1)^{2m}+2^{2m}\hat{\lambda}}{(l^{2m}+(k+1)^{2m})(l^{2m}+k^{2m})+2^{2m}\hat{\lambda}(l^{2m}+k^{2m})+\hat{\lambda}^2}\leq \frac{2^{2m}}{\hat{\lambda}},
\end{equation}
where $k=N-2j$. The last inequality reduces to
\begin{gather*}
  \frac{l^{2m}+2^{2m-1}(k^{2m}+1)+2^{2m}\hat{\lambda}}{(l^{2m}+k^{2m})^2+2^{2m}\hat{\lambda}(l^{2m}+k^{2m})+\hat{\lambda}^2}\leq \frac{2^{2m}}{\hat{\lambda}},
\end{gather*}
after grouping the terms in this inequality it is easy to see that it is satisfied for all $l\geq 1$, and $k\geq 0$.

Obviously $\hat{\lambda}^{-1}\geq b_1$ holds. Analogically as above, assuming $\hat{\lambda}^{-1}\geq b_{j}$ we verify that $\hat{\lambda}^{-1}\geq b_{j+1}$ holds (it can also be verified that $f^\prime(b)$ is strictly positive, and it is enough to verify the inequality setting $b_j = \hat{\lambda}^{-1}$).
\begin{equation}
  \label{eq:seriesBbound}
  \frac{d_{2j+1}+b_{j-1}\hat{\lambda}^2}{d_{2j+2}d_{2j+1}+b_{j-1}d_{2j+2}\hat{\lambda}^2 + \hat{\lambda}^2}\leq\frac{l^{2m}+(2j)^{2m}+\hat{\lambda}}{(l^{2m}+(2j+1)^{2m})(l^{2m}+(2j)^{2m})+\hat{\lambda}(l^{2m}+(2j+1)^{2m}) + \hat{\lambda}^2}\leq \frac{1}{\hat{\lambda}},
\end{equation}
The last inequality holds due to following inequality, which is clearly satisfied
\begin{gather*}
  \frac{a+\hat{\lambda}}{a^2+a\hat{\lambda}+\hat{\lambda}^2}\leq \frac{1}{\hat{\lambda}},
\end{gather*}
where $a=l^{2m}+(2j)^{2m}$.

\qed

\providecommand{\llambda}{\hat{\lambda}}

\begin{lemma}
  \label{lem:bounds_diagonal}
  Let $l>1$. All elements in $2\times 2$ diagonal blocks of $\operS^{-1}$ are estimated uniformly. Precisely, the following inequalities hold for $j=0,\dots,N/2-1$
  \begin{align*}
    |(\operS^{-1})_{2j+1,2j+1}|\leq \seriesAconst / \hat{\lambda},\qquad|(\operS^{-1})_{2j+1,2j+2}|\leq 1 / \hat{\lambda},\\
    |(\operS^{-1})_{2j+2,2j+1}|\leq 1 / \hat{\lambda},\qquad|(\operS^{-1})_{2j+2,2j+2}|\leq 1 / \hat{\lambda}.
  \end{align*}
\end{lemma}
\paragraph{Proof}
Here we assume that $l$ is fixed, and we drop the superscript in the notation of $\matrixS^l$, $\seriesA{k}$, $\seriesB{k}$, and $\seriesD{k}$, we use simply $\matrixS$, $a_k$, $b_k$, and $d_k$ respectively. First, all matrices considered are invertible, which is obvious by calculating the determinant of tridiagonal matrices.

\newcommand{\W}{\mbox{Inv}}
\newcommand{\V}{\overline{\mbox{Inv}}}
We use the 
notation $\W^{(D)}$ to denote the $D\times D$ dimensional upper-left corner block of 
$\operS^{-1}$. Analogously we use the notation $\V^{(D)}$ to denote the $D\times D$ dimensional lower-right corner block of $\operS^{-1}$. 
We are going to use the following convention for block decomposition of $\operS^{-1} = \W^{(N)}$.
\begin{subequations}
\begin{align}
  \operS^{-1}=\left[\begin{array}{cc}\W^{(N)}_{11}&\W^{(N)}_{12}\\ \W^{(N)}_{21}& \W^{(N)}_{22}\end{array}\right]&:=\left[\begin{array}{cc}\matrixS(d_1,\dots,d_{N-2})&-A\\A^T&\matrixS(d_{N-1},d_N)\end{array}\right]^{-1},\label{eq:Wdecomposition}\\
  \operS^{-1}=\left[\begin{array}{cc}\V^{(N)}_{11}&\V^{(N)}_{12}\\ \V^{(N)}_{21}&\V^{(N)}_{22}\end{array}\right]&:=\left[\begin{array}{cc}\matrixS(d_1,d_{2})&-{A^\prime}^T\\A^\prime&\matrixS(d_3,\dots,d_{N-1},d_N)\end{array}\right]^{-1},\label{eq:Vdecomposition}
\end{align}
\end{subequations}
where
\begin{equation*}
  A=\left[\begin{array}{cc}0&0\\\vdots&\vdots\\0&0\\\llambda&0\end{array}\right],\qquad
  A^\prime=\left[\begin{array}{cc}0&\llambda\\\vdots&\vdots\\0&0\\0&0\end{array}\right].
\end{equation*}
We will call $\W^{(N)}_{j,k}$, $\V^{(N)}_{j,k}$ the inverse blocks. In the remainder of the proof we will compute recursively $\W^{(N-2)}$, $\V^{(N-2)}$ ... $\W^{(2)}$, $\V^{(2)}$.

The explicit formulas for the inverse blocks are obtained from the following system of equations (simplifying the notation by
dropping the brackets with parameters, i.e. $\matrixS_I=\matrixS(d_1,\dots,d_{N-2}),\ \matrixS_{II}=\matrixS(d_{N-1},d_N)$ etc.)
\begin{subequations}
\label{eq:inverse_blocks}
\begin{align}  
  \matrixS_I\W^{(N)}_{11} - A\W^{(N)}_{21}=I,&\quad&\matrixS_I\W^{(N)}_{12} - A\W^{(N)}_{22}=0,\\A^T\W^{(N)}_{11}+\matrixS_{II}\W^{(N)}_{21}=0,&\quad&A^T\W^{(N)}_{12}+\matrixS_{II}\W^{(N)}_{22}=I.
\end{align}
\end{subequations}
When the equations for diagonal blocks are decoupled we obtain
\begin{align*}
  \W^{(N)}_{11}&=\left[\matrixS(d_1,\dots,d_{N-2})+A\matrixS(d_{N-1},d_N)^{-1}A^T\right]^{-1},\\ \W^{(N)}_{22}&=\left[\matrixS(d_{N-1},d_N)+A^T\matrixS(d_1,\dots,d_{N-2})^{-1}A\right]^{-1}.
\end{align*}
where
\begin{align*}
  A\matrixS(d_{N-1},d_N)^{-1}A^T&=\left[\begin{array}{cccc}0&\dots&0&0\\\vdots&\ &\vdots&\vdots\\0&\dots&0&0\\0&\dots&0&\llambda^2\matrixS(d_{N-1},d_N)^{-1}_{11}\end{array}\right],\\
  A^T\matrixS(d_1,\dots,d_{N-2})^{-1}A&=\left[\begin{array}{cc}\llambda^2\matrixS(d_1,\dots,d_{N-2})^{-1}_{NN}&0\\0&0\end{array}\right].
\end{align*}
Now, we state the crucial observation -- the inverse diagonal blocks $W_{11}$ and $W_{22}$ are inverses of tridiagonal matrices, i.e.
\begin{subequations}
  \label{eq:Wdecoupling}
  \begin{align}
    \W^{(N)}_{11}&=\matrixS(d_1,\dots,d_{N-3},\,d_{N-2}+\llambda^2a_1)^{-1} = \W^{(N-2)},\\
    \W^{(N)}_{22}&=\matrixS(d_{N-1}+\llambda^2\matrixS(d_1,\dots,d_{N-2})^{-1}_{NN},\,d_N)^{-1},
  \end{align}
\end{subequations}
where $a_1=\matrixS(d_{N-1},d_N)^{-1}_{11}=\frac{d_N}{d_{N-1}d_N+\llambda^2}$. The same holds for the diagonal inverse blocks $\V^{(N)}_{11}$ and $\V^{(N)}_{22}$ by symmetric calculations, i.e.
\begin{subequations}
  \label{eq:Vdecoupling}
  \begin{align}
    \V^{(N)}_{11}&=\left[\matrixS(d_1,d_2)+{A^\prime}^T\matrixS(d_3,\dots,d_N)^{-1}{A^\prime}\right]^{-1}=\matrixS\left(d_1,\,d_2+\llambda^2\matrixS(d_3,\dots,d_N)^{-1}_{11}\right)^{-1},\\
    \V^{(N)}_{22}&=\left[\matrixS(d_3,\dots,d_N)+A^\prime\matrixS(d_1,d_2)^{-1}{A^\prime}^T\right]^{-1}=\matrixS(d_3+\llambda^2b_1,\,d_4,\dots,d_N)^{-1},
  \end{align}
\end{subequations}
where $b_1=\matrixS(d_1,d_2)^{-1}_{22}=\frac{d_1}{d_{1}d_2+\llambda^2}$.

Observe that the decoupling of the diagonal blocks described above can be iterated, and the matrix $\W^{(N-2)}=\matrixS(d_1,\dots,d_{N-3},d_{N-2}+\llambda^2a_1)^{-1}$ is further decomposed
\begin{equation*}
  \W^{(N-2)} = \left[\begin{array}{cc}\W^{(N-2)}_{11}&\W^{(N-2)}_{12}\\ \W^{(N-2)}_{21}&\W^{(N-2)}_{22}\end{array}\right] = \left[\begin{array}{cc}\matrixS(d_1,\dots,d_{N-4})&-A\\A^T&\matrixS(d_{N-3},d_{N-2}+\llambda^2a_1)\end{array}\right]^{-1},
\end{equation*}
thus  we write the formula for the inverse diagonal block $\W^{(N-2)}_{11} = \W^{(N-4)}$
\begin{multline*}
  \W^{(N-2)}_{11}= \W^{(N-4)}=\left[\matrixS(d_1,\dots,d_{N-4})+A\matrixS(d_{N-3},d_{N-2}+\llambda^2a_1)^{-1}A^T\right]^{-1}\\=\matrixS(d_1,\dots,d_{N-5},\,d_{N-4}+\llambda^2a_2)^{-1},
\end{multline*}
where
\begin{equation*}
  a_2=\matrixS(d_{N-3},d_{N-2}+\llambda^2a_1)^{-1}_{11}=\frac{d_{N-2}+a_1\llambda^2}{d_{N-3}d_{N-2}+a_1d_{N-3}\llambda^2+\llambda^2}.
\end{equation*}

From repeating $j$ times the procedure of taking the upper-left inverse diagonal block and decompose it further like in \eqref{eq:Wdecomposition}, we obtain the explicit formula for the $N-2j$ dimensional upper-left diagonal block of $\operS^{-1}$
  \begin{align}
    \label{eq:upperLeftBlock}
    \W^{(N-2j)}&=\left[\matrixS(d_1,\dots,d_{N-2j})+A\matrixS(d_{N-2j+1},d_{N-2j+2}+\llambda^2a_{j-1})^{-1}A^T\right]^{-1}\\&=\matrixS(d_1,\dots,d_{N-2j-1},\,d_{N-2j}+\llambda^2a_{j})^{-1},
  \end{align}
where
\begin{equation*}
  a_j=\matrixS(d_{N-2j+1},d_{N-2j+2}+\llambda^2a_{j-1})^{-1}_{11}=\frac{d_{N-2j+2}+a_{j-1}\llambda^2}{d_{N-2j+1}d_{N-2j+2}+a_{j-1}d_{N-2j+1}\llambda^2+\llambda^2}.
\end{equation*}
Performing iteratively $j$ times the symmetric procedure to the one described above (performing decomposition like in \eqref{eq:Vdecomposition}), we obtain the explicit formula for $N-2j$ dimensional lower right diagonal inverse block
  \begin{align}
    \V^{N-2j}&=\left[\matrixS(d_{2j+1},\dots,d_{N})+A^\prime\matrixS(d_{2j-1}+\llambda^2b_{j-1},d_{2j})^{-1}{A^\prime}^T\right]^{-1}\\&=\matrixS(d_{2j+1}+\llambda^2b_{j},\,d_{2j+2},\dots,d_{N})^{-1},
  \end{align}
where
\begin{equation*}
  b_j=\matrixS(d_{2j-1}+\llambda^2b_{j-1},d_{2j})^{-1}_{22}=\frac{d_{2j-1}+\llambda^2b_{j-1}}{d_{2j-1}d_{2j}+b_{j-1}d_{2j}\llambda^2+\llambda^2}.
\end{equation*}
Note that the recursive series $\{a_j\}$, $\{b_j\}$ are generated from the procedures described above.
%
%
%
%
Using above results, we may now derive an explicit formulas for the $2\times 2$ diagonal blocks of $\operS^{-1}$. Let us present an example how it is done.
Observe that from \eqref{eq:upperLeftBlock} we have that the $N-2j$ dimensional upper-left block of $\operS^{-1}$ is $\W^{(N-2j)}=\matrixS(d_1,\dots,d_{N-2j-1},d_{N-2j}+\llambda^2a_{j})^{-1}$.

Then, for $j=(N-2)/2$ we are left with $\W^{(2)}$, whereas if $j < (N-2)/2$ we apply $j$ times to $\W^{(N-2j)}$ the procedure of taking the lower right diagonal block, and decomposing like in \eqref{eq:Vdecomposition}, and we get that the $j$-th (counting from the bottom) $2\times 2$ diagonal block of $\operS^{-1}$ equals to
\begin{equation*}
  \matrixS(d_{N-2j-1}+\llambda^2b_{k}, d_{N-2j}+\llambda^2a_{j})^{-1}=\left[\begin{array}{cc}d_{N-2j-1}+\llambda^2b_{k}&-\llambda\\\llambda&d_{N-2j}+\llambda^2a_{j}\end{array}\right]^{-1},\text{ where }k=(N-2j-2)/2.
\end{equation*}

Let us denote $D=\llambda^2+d_{N-2j-1}d_{N-2j}+d_{N-2j}b_k\llambda^2+d_{N-2j-1}a_j\llambda^2+\llambda^4a_jb_k$, we have
\begin{align*}
  \matrixS(d_{N-2j-1}+\llambda^2b_{k}, d_{N-2j}+\llambda^2a_{j})^{-1}_{11}=&\frac{d_{N-2j}+\llambda^2a_j}{D},\\
  \matrixS(d_{N-2j-1}+\llambda^2b_{k}, d_{N-2j}+\llambda^2a_{j})^{-1}_{22}=&\frac{d_{N-2j-1}+\llambda^2b_{k}}{D},\\
  \begin{array}{c}\matrixS(d_{N-2j-1}+\llambda^2b_{k}, d_{N-2j}+\llambda^2a_{j})^{-1}_{12}\\\matrixS(d_{N-2j-1}+\llambda^2b_{k}, d_{N-2j}+\llambda^2a_{j})^{-1}_{21}\end{array}=&\frac{\pm\llambda}{D}.
\end{align*}
We have that for $f(a_j,b_k)=\frac{d_{N-2j}+\llambda^2a_j}{D}=\frac{N}{D}$ the partial derivatives equal to
\begin{equation*}
  \frac{\partial f(a_j,b_k)}{\partial a_j}=\frac{\llambda^4}{D^2}>0,\qquad\frac{\partial f(a_j,b_k)}{\partial b_k}=\frac{-\llambda^2N^2}{D^2}<0.
\end{equation*}
Hence, to bound $\matrixS(d_{N-2j-1}+\llambda^2b_{k}, d_{N-2j}+\llambda^2a_{j})^{-1}_{11}$ we use the upper end of the bound for $a_j$ from Lemma~\ref{lem:sequences}, i.e. we set $a_j=\seriesAconst\llambda^{-1}$, and we use the lower end of the bound for $b_k$, i.e. we set $b_k=0$.
We are left with bounding $\frac{d_{N-2j}+\seriesAconst \llambda}{\llambda^2+d_{N-2j-1}d_{N-2j}+d_{N-2j-1}\seriesAconst \llambda}$,
which was already showed in \eqref{eq:seriesAbound} to be bounded by $\seriesAconst\llambda^{-1}$.

To bound $\matrixS(d_{N-2j-1}+\llambda^2b_{k}, d_{N-2j}+\llambda^2a_{j})^{-1}_{22}$, analogously as above, we set
$a_j=0$, and $b_k=\llambda^{-1}$, and we are left with bounding $\frac{d_{N-2j-1}+\llambda}{\llambda^2+d_{N-2j-1}d_{N-2j}+d_{N-2j}\llambda}$, which was already showed in \eqref{eq:seriesBbound} to be bounded by $\llambda^{-1}$.
To bound the remaining two elements, i.e. $\matrixS(d_{N-2j-1}+\llambda^2b_{k}, d_{N-2j}+\llambda^2a_{j})^{-1}_{12}$, $\matrixS(d_{N-2j-1}+\llambda^2b_{k}, d_{N-2j}+\llambda^2a_{j})^{-1}_{21}$ we set $a_j=0$, $b_k=0$, and we obtain the claimed bounds immediately.

\qed

\begin{lemma}
\label{lem:uniformBound}
  Let $l>0$. The following uniform bound hold
\begin{subequations}
    \label{lambdaBounds}
  \begin{equation}
    \begin{array}{ll}
    |(\operS^{-1})_{2j+1,2k+1}|\leq \seriesAconst / \hat{\lambda}&|(\operS^{-1})_{2j+1,2k+2}|\leq 1 / \hat{\lambda}\\
    |(\operS^{-1})_{2j+2,2k+1}|\leq 1 / \hat{\lambda}&|(\operS^{-1})_{2j+2,2k+2}|\leq 1 / \hat{\lambda}
    \end{array}
    \text{, for }k,j=0,\dots,[(N-1)/2].
  \end{equation}
\end{subequations}
\end{lemma}
\paragraph{Proof}
Here we assume that $l$ is fixed, and we drop the superscript in the notation of $\matrixS^l$, $\seriesA{k}$, $\seriesB{k}$, and $\seriesD{k}$, we use simply $\matrixS$, $a_k$, $b_k$, and $d_k$ respectively. We use the same 
notation as in Lemma~\ref{lem:bounds_diagonal}, i.e. we use $\W^{(D)}$ to denote the $D\times D$ dimensional upper-left corner block of 
$\operS^{-1}$. Analogously we use the notation $\V^{(D)}$ to denote the $D\times D$ dimensional lower-right corner block of $\operS^{-1}$. 

First, for the sake of presentation, let us prove that the claimed bounds are true for the $4\times 4$ upper left corner submatrix of
$\operS^{-1}$, i.e. $\W^{(4)}$, the general result will follow
\begin{equation}
\label{eqw4}
  \W^{(4)}=\matrixS(d_1,d_{2},d_{3},d_4+\llambda^2a_{N/2-2})^{-1}=\left[\begin{array}{cc}\W^{(4)}_{11}&\W^{(4)}_{12}\\ \W^{(4)}_{21}&\W^{(4)}_{22}\end{array}\right]:=\left[\begin{array}{cc}\matrixS(d_1,d_{2})&-A\\A^T&\matrixS(d_{3},d_4+\llambda^2a_{N/2-2})\end{array}\right]^{-1}.
\end{equation}
From the equations for inverse blocks \eqref{eq:inverse_blocks} it follows that the block beyond diagonal satisfies 
\begin{equation*}
  \W^{(4)}_{21}=\matrixS(d_{3},d_4+\llambda^2a_{N/2-2})^{-1}A^T\W^{(4)}_{11}.
\end{equation*}
From Lemma~\ref{lem:sequences} follows that $\seriesAconst\llambda^{-1}\geq a_{N/2-2}\geq0$, hence the bounds for all elements of $\matrixS(d_{3},d_4+\llambda^2a_{N/2-2})^{-1}$ are the same as those derived in Lemma~\ref{lem:bounds_diagonal}. Observe that
\begin{multline*}
  \W^{(4)}_{21}=\matrixS(d_{3},d_4+\llambda^2a_{N/2-2})^{-1}A^T\W^{(4)}_{11}=\\\left[\begin{array}{cc}\matrixS(d_{3},d_4+\llambda^2a_{N/2-2})^{-1}_{11}\llambda(\W^{(4)}_{11})_{21}&\matrixS(d_{3},d_4+\llambda^2a_{N/2-2})^{-1}_{11}\llambda(\W^{(4)}_{11})_{22}\\\matrixS(d_{3},d_4+\llambda^2a_{N/2-2})^{-1}_{21}\llambda(\W^{(4)}_{11})_{21}&
\matrixS(d_{3},d_4+\llambda^2a_{N/2-2})^{-1}_{21}\llambda(\W^{(4)}_{11})_{22}\end{array}\right],
\end{multline*}
as we have from Lemma~\ref{lem:bounds_diagonal} the bounds $|\llambda(\W^{(4)}_{11})_{21}|,\,|\llambda(\W^{(4)}_{11})_{22}|\leq 1$, and
\\$|\matrixS(d_{3},d_4+\llambda^2a_{N/2-2})^{-1}_{11}|\leq\seriesAconst\llambda^{-1},\ |\matrixS(d_{3},d_4+\llambda^2a_{N/2-2})^{-1}_{21}|\leq\llambda^{-1}$ . Elements from the block $\W^{(4)}_{21}$ clearly satisfy the following bounds
\begin{equation*}
  |(\W^{(4)}_{21})_{11}|,\,|(\W^{(4)}_{21})_{12}|\leq \seriesAconst\llambda^{-1},\qquad |(\W^{(4)}_{21})_{21}|,\,|(\W^{(4)}_{21})_{22}|\leq \llambda^{-1}.
\end{equation*}
The block $\W^{(4)}_{12}$ satisfies symmetric bounds by a symmetric argument. Observe that the bounds for the diagonal blocks $\W^{(4)}_{11}$, $\W^{(4)}_{22}$ were
derived in the previous lemma, hence at this point we have bounded uniformly all elements in $\W^{(4)}=\matrixS(d_1,d_{2},d_{3},d_4+\llambda^2a_{N/2-2})^{-1}$.
Observe that in order to derive the bounds for the off-diagonal blocks, we used only the bound for the last row of $\W^{(4)}$ ($(\W^{(4)}_{11})_{21}=\W^{(4)}_{41},\ (\W^{(4)}_{11})_{22}=\W^{(4)}_{42},\ (\W^{(4)}_{22})_{21}=\W^{(4)}_{43},\ (\W^{(4)}_{22})_{22}=\W^{(4)}_{44}$ see \eqref{eqw4}). From the bounds established so far all elements in the last row of $\W^{(4)}$ satisfy $|\W^{(4)}_{4j}|\leq\llambda^{-1}$.
It is easy to see that, if we now consider $\W^{(6)}=\matrixS(d_1,d_2,d_{3},d_4,d_5,d_6+\llambda^2a_{N/2-3})^{-1}$,
by a similar argument for $j=1,\dots,6$ we obtain the bounds
\begin{equation*}
\begin{array}{ll}
\left|\W^{(6)}_{6j}\right|\leq\llambda^{-1},\\
\left|\W^{(6)}_{5j}\right|\leq \seriesAconst\llambda^{-1}.
\end{array}
\end{equation*}


Finally, from the presentation above follows that assuming that absolute value of all of the elements in the last row of the inverse block $\W^{(2k)}=\matrixS(d_1,\dots,d_{2k-1},d_{2k}+\llambda^2a_{N/2-k})^{-1}$ are bounded by $\llambda^{-1}$, and the rest by $\seriesAconst\llambda^{-1}$,
the same bounds for the larger inverse block $\W^{(2k+2)}=\matrixS(d_1,\dots,d_{2k+1},d_{2k+2}+\llambda^2a_{N/2-(k+1)})^{-1}$ will follow,
thus, we showed that the bounds \eqref{lambdaBounds} are propagated for the whole $\operS^{-1}=\W^{(N)}$.

\qed

\begin{lemma}
  \label{lem:columnEstimate}
  Let $N>0$, $m>1$, $l=1,\dots,N$. There exist $C(m)>0$ (independent of $\lambda$ and $N$), such that for
  $\lambda>1$ the following bounds hold
  \begin{align*}
    \sum_{j=1,\dots,N}{\left|(\operS^{-1})_{ij}\right|}&\leq C(m)\cdot \hat{\lambda}^{-1+1/2m},\\
    \sum_{j=1,\dots,N}{\left|(\operS^{-1})_{ji}\right|}&\leq C(m)\cdot \hat{\lambda}^{-1+1/2m},
  \end{align*}
  for all $i=1,\dots,N$.
\end{lemma}

\paragraph{Proof}
Here we assume that $l$ is fixed, and we drop the superscript in the notation of $\matrixS^l$, $\seriesA{k}$, $\seriesB{k}$, and $\seriesD{k}$, we use simply $\matrixS$, $a_k$, $b_k$, and $d_k$ respectively.

For the sake of clarification let us restrict our attention to the first column of $\operS^{-1}$.

From Lemma~\ref{lem:bounds_diagonal} it follows that $n=[c\llambda^{1/2m}]$
dimensional upper left corner submatrix of $\operS^{-1}$ is equal to $ \W^{(n)} = \matrixS(d_1,\dots, d_{n-1}, d_n+\llambda^2a_{(N-n)/2})^{-1} $, and the absolute values of elements
in this matrix are uniformly bounded by $\seriesAconst\llambda^{-1}$, thus the straightforward estimate for the first part of the sum is
\begin{equation*}
  \sum_{i=1,\dots,n}{\left|(\operS^{-1})_{i1}\right|}\leq n\seriesAconst\llambda^{-1} =  c\llambda^{1/2m}\cdot\seriesAconst\llambda^{-1}= {\seriesAconst c}{\llambda^{-1+1/2m}}.
\end{equation*}
Next, we are going to show that the terms in remainder 
\[
\sum_{i=n+1,\dots,N}{\left|(\operS^{-1})_{i1}\right|}=\sum_{i=n+1,\dots,N}{\left|\W^{(N)}_{i1}\right|}
\]
obey a geometric decay rate, and can be bounded uniformly with respect to the dimension.

As in the previous lemmas we take the block decomposition of 
\[
\operS^{-1} = \W^{(N)}\text{, i.e.}
\]
\begin{equation*}
  \W^{(N)}=\left[\begin{array}{cc}\W^{(N)}_{11}&\W^{(N)}_{12}\\ \W^{(N)}_{21}&\W^{(N)}_{22}\end{array}\right]=\left[\begin{array}{cc}\matrixS(d_1,\dots,d_{N-2})&-A\\A^T&\matrixS(d_{N-1},d_N)\end{array}\right]^{-1}.
\end{equation*}
From \eqref{eq:inverse_blocks} it follows that $\W^{(N)}_{21}$ can be expressed in terms of $\W^{(N)}_{11}$, namely, for the first column of $\W^{(N)}_{21}$ 
the identities are
\begin{equation*}
  \W^{(N)}_{N-1,1}=\matrixS(d_{N-1},d_N)^{-1}_{11}\llambda\W^{(N)}_{N-2,1},\quad \W^{(N)}_{N,1}=\matrixS(d_{N-1},d_N)^{-1}_{21}\llambda\W^{(N)}_{N-2,1}.
\end{equation*}

Analogously, for $\W^{(N)}_{N-3,1},\ \W^{(N)}_{N-2,1}$ we have
\begin{equation*}
  \W^{(N)}_{N-3,1}=\matrixS(d_{N-3},d_{N-2}+\llambda^2a_1)^{-1}_{11}\llambda\W^{(N)}_{N-4,1},\quad \W^{(N)}_{N-2,1}=\matrixS(d_{N-3},d_{N-2}+\llambda^2a_1)^{-1}_{21}\llambda\W^{(N)}_{N-4,1}.
\end{equation*}
From repeating this argument we obtain
\begin{align*}
  &\W^{(N)}_{N-5,1}=\matrixS(d_{N-5},d_{N-4}+\llambda^2a_2)^{-1}_{11}\llambda\W^{(N)}_{N-6,1},\quad \W^{(N)}_{N-4,1}=\matrixS(d_{N-5},d_{N-4}+\llambda^2a_2)^{-1}_{21}\llambda\W^{(N)}_{N-6,1}, \\
&\dots,\\
  &\W^{(N)}_{n+1,1}=\matrixS(d_{n+1},d_{n+2}+\llambda^2a_{(N-n)/2-1})^{-1}_{11}\llambda\W^{(N)}_{n1},\quad \W^{(N)}_{n+2,1}=\matrixS(d_{n+1},d_{n+2}+\llambda^2a_{(N-n)/2-1})^{-1}_{21}\llambda\W^{(N)}_{n1}.
\end{align*}
This is a recursive series, all elements can be expressed in terms of $\W^{(N)}$. Therefore
\begin{equation*}
  \sum_{i=n+1,\dots,N}{\left|\W^{(N)}_{i1}\right|}\leq\left|\W^{(N)}_{n1}\right|\sum_{j=n+1}^{N}{|c_j|},
\end{equation*}
where
\begin{align*}
  c_{n+1}&:=\matrixS(d_{n+1},d_{n+2}+\llambda^2a_{(N-n)/2-1})^{-1}_{11}\llambda,\\
  c_{n+2}&:=\matrixS(d_{n+1},d_{n+2}+\llambda^2a_{(N-n)/2-1})^{-1}_{21}\llambda,\\
  c_{n+3}&:=\matrixS(d_{n+3},d_{n+4}+\llambda^2a_{(N-n)/2-2})^{-1}_{11} \llambda\cdot\matrixS(d_{n+1},d_{n+2}+\llambda^2a_{(N-n)/2-1})^{-1}_{21}\llambda ,\\
  c_{n+4}&:=\matrixS(d_{n+3},d_{n+4}+\llambda^2a_{(N-n)/2-2})^{-1}_{21} \llambda\cdot\matrixS(d_{n+1},d_{n+2}+\llambda^2a_{(N-n)/2-1})^{-1}_{21}\llambda,\\
  &\dots,\\
  c_{N-1}&:=\matrixS(d_{N-1},d_N)^{-1}_{11}\llambda\cdot\matrixS(d_{N-3},d_{N-2}+\llambda^2a_1)^{-1}_{21}\llambda\cdots \matrixS(d_{n+1},d_{n+2}+\llambda^2a_{(N-n)/2-1})^{-1}_{21}\llambda,\\
  c_{N}&:=\matrixS(d_{N-1},d_N)^{-1}_{21}\llambda\cdot\matrixS(d_{N-3},d_{N-2}+\llambda^2a_1)^{-1}_{21}\llambda\cdots \matrixS(d_{n+1},d_{n+2}+\llambda^2a_{(N-n)/2-1})^{-1}_{21}\llambda.
\end{align*}

From Lemma~\ref{lem:sequences} it follows that $\seriesAconst\llambda^{-1}\geq a_{j}\geq 0$, and $f(a)=\frac{d_{j+1}+\llambda^2 a}{\llambda^2+d_{j}d_{j+1}+d_{j}\llambda^2 a }$ is strictly positive for all $j\geq 1$ (as the derivative is positive, compare \eqref{eq:strictlyIncreasing}). Recall that $d_j=l^{2m}+(j-1)^{2m}$, and we have the obvious inequality $\llambda^2+d_{j}d_{j+1}+d_{j}\llambda^2 a<\llambda^2+d_{j+1}d_{j+2}+d_{j+1}\llambda^2 a$ for $a>0$. Therefore the following inequalities are satisfied
\begin{align*}
  \left|\matrixS(d_{n+j},d_{n+j+1}+\llambda^2a_{(N-(n+j+1))/2})^{-1}_{11}\right|\leq\left|\matrixS(d_{n+1},d_{n+2}+\llambda^2a_{(N-n)/2-1})^{-1}_{11}\right|&\leq\frac{d_{n+2}+2^{2m}\llambda }{\llambda^2+d_{n+1}d_{n+2}+2^{2m}\llambda d_{n+1}},\\
  \left|\matrixS(d_{n+j},d_{n+j+1}+\llambda^2a_{(N-(n+j+1))/2})^{-1}_{21}\right|\leq\left|\matrixS(d_{n+1},d_{n+2}+\llambda^2a_{(N-n)/2-1})^{-1}_{21}\right|&\leq\frac{\llambda}{\llambda^2+d_{n+1}d_{n+2}}.
\end{align*}

Now to show the claim about the geometric decay, we take $n>[2^{(2m+1)/2m}\llambda^{1/2m}]$
\begin{multline*}
  \frac{\llambda}{\llambda^2+d_{n+1}d_{n+2}}\leq\frac{d_{n+2}+2^{2m}\llambda}{\llambda^2+d_{n+1}d_{n+2}+2^{2m}\llambda d_{n+1}}=\frac{l^{2m}+(n+1)^{2m}+2^{2m}\llambda}{\llambda^2+(l^{2m}+n^{2m})(l^{2m}+(n+1)^{2m})+2^{2m}\llambda(l^{2m}+n^{2m})}
\end{multline*}
Similar argument to the one used in \eqref{eq:abestimate} shows that
\begin{multline*}
  \frac{l^{2m}+(n+1)^{2m}+2^{2m}\llambda}{\llambda^2+(l^{2m}+n^{2m})(l^{2m}+(n+1)^{2m})+2^{2m}\llambda(l^{2m}+n^{2m})}\leq
  \frac{\seriesAconst a+\seriesAconst\llambda}{a^2+\seriesAconst a\llambda+\llambda^2}\leq
  \frac{2^{4m+1}\llambda+2^{2m}\llambda}{2^{4m+2}\llambda^2+2^{4m+1}\llambda^2+\llambda^2}<\frac{1}{2\llambda}.
\end{multline*}

We thus demonstrated that for any $N$ we have
\begin{equation}
  \label{eq:geometricRegime}
  \sum_{i=n+1,\dots,N}{\left|(\operS^{-1})_{i1}\right|}=\sum_{i=n+1,\dots,N}{\left|\W^{(N)}_{i1}\right|}<\left|\W^{(N)}_{n1}\right|2\sum_{i=1}^{\infty}{\left(\frac{1}{2}\right)^i}=2\left|\W^{(N)}_{n1}\right|\leq\frac{2^{2m+1}}{\llambda}.
\end{equation}
Now taking $C(m)>2^{(2m+1)/2m}+2^{2m+1}$ we obtain the claim.

To conclude, observe that the bound holds for the first row, as the matrices $\operS$ and $\operS^{-1}$ commute ($\operS=D+A$, where $D$ is a diagonal, $A$ is a skew-symmetric matrix).
The bound is true for any other column/row, to see this note that for each column there are at most $n=[c\llambda^{1/2m}]$ elements
beyond the geometric decay regime, therefore the bound is true for any column of $\operS^{-1}$.

\qed

\subsection{Proof of Theorem~\ref{lem:invMatrixEstimates}}
Using lemmas presented in this section we prove the main result with inverse matrix bounds

\begin{reptheorem}{lem:invMatrixEstimates}
  Let $l=1,\dots,N$. Let $\operSS{l}$ be the matrix given by \eqref{eq:Llambdal}, $\operSS{B}$ be
the matrix $(P_BL_\lambda)^{-1}$.

  The following estimates hold for the matrices $\operSS{l}^{-1}$ (diagonal submatrices of $\operSS{B}^{-1}$).
  \begin{align*}
     \|\operSS{l}^{-1}\|_{l^1\to l^1}&\leq C_1(m)\left(\frac{l\lambda}{2}\right)^{-1+1/2m},\\
     \|\operSS{l}^{-1}\|_{l^1\to l^\infty}&\leq\seriesAconst\left(\frac{l\lambda}{2}\right)^{-1},\\
     \|\operSS{l}^{-1}\|_{l^1\to l^1_1}&\leq C_2(m)\left(\frac{l\lambda}{2}\right)^{-1+1/m}.
   \end{align*}

  The following estimates hold for the matrix $\operSS{B}^{-1}$
   \begin{align*}
     \|\operSS{B}^{-1}\|_{l^1\to l^1}&\leq C_1(m)\left(\frac{\lambda}{2}\right)^{-1+1/2m},\\
     \|\operSS{B}^{-1}\|_{l^1\to l^\infty}&\leq\seriesAconst\left(\frac{\lambda}{2}\right)^{-1},\\
     \|\operSS{B}^{-1}\|_{l^1\to l^1_1}&\leq C_2(m)\left(\frac{\lambda}{2}\right)^{-1+1/m}.
   \end{align*}
\end{reptheorem}

\paragraph{Proof} The uniform $l^1$ estimate for each column of $\operSS{B}^{-1}$ follows directly from Lemma~\ref{lem:columnEstimate}.
The uniform $l^\infty$ estimate for each column of $\operSS{B}^{-1}$ follows directly from Lemma~\ref{lem:uniformBound}.

In order to estimate uniformly the gradient norm of $\operSS{B}^{-1}$ we are going to consider two cases separately.

Let $c>2^{(2m+1)/2m}$, $\alpha=\left[c\left(\frac{\lambda}{2}\right)^{1/(2m-1)}\right]$, where $[\cdot]$ is the integer part. 
Let us demonstrate the result for the first column of $\operSS{B}^{-1}$.
\paragraph{Case I} For $l\leq\alpha$ we split the sum
\begin{equation*}
  \sum_{k=1}^N{(l+k)\left|\left(\operS^{-1}\right)_{k,1}\right|}=\sum_{k=1}^n{\left|\left(\operS^{-1}\right)_{k,1}\right|}+
  \sum_{k=n+1}^N{(l+k)\left|\left(\operS^{-1}\right)_{k,1}\right|},
\end{equation*}
where $n=l^{1/2m}\left(\frac{\lambda}{2}\right)^{-1/[2m(2m-1)]}\alpha=\left[c\left(\frac{l\lambda}{2}\right)^{1/2m}\right]$ (this particular choice is due to technical reasons).
The finite part of the sum above can be estimated

\begin{equation*}
  \sum_{k=1}^n{\left|\left(\operS^{-1}\right)_{k,1}\right|}\leq\max_{k=1,\dots,N}\left|\left(\operS^{-1}\right)_{k,1}\right|\sum_{k=1}^n{l+k}
  \leq 2^{2m}\left(\frac{l\lambda}{2}\right)^{-1}\left(ln+n^2\right)\leq c_1\left(\frac{l\lambda}{2}\right)^{-1}\left(\frac{l\lambda}{2}\right)^{1/m}=c_1\left(\frac{l\lambda}{2}\right)^{-1+1/m},
\end{equation*}
where  we estimated $ln=l^{(2m-1)/2m}l^{1/m}\left(\frac{\lambda}{2}\right)^{1/2m}\leq\alpha^{(2m-1)/2m}l^{1/m}\left(\frac{\lambda}{2}\right)^{1/2m}\leq \tilde{c} \left(\frac{l\lambda}{2}\right)^{1/m}$.
As $d^l_n>2^{2m+1}\frac{l\lambda}{2}$ from the proof of Lemma~\ref{lem:columnEstimate} it follows that the remaining part of the sum
is within the geometric decay regime, therefore we can estimate like in \eqref{eq:geometricRegime}
\begin{equation*}
  \sum_{k=n+1}^N{(l+k)\left|\left(\operS^{-1}\right)_{k,1}\right|}\leq
\left|\left(\operS^{-1}\right)_{n,1}\right|2\left(l\sum_{k=1}^\infty{\frac{1}{2^k}+\sum_{k=1}^\infty{\frac{k}{2^k}}}\right)\leq
c_2\left|\left(\operS^{-1}\right)_{n,1}\right|(l+1)\leq C_2\left(\frac{\lambda}{2}\right)^{-1},
\end{equation*}
in the last inequality we used the estimate from Lemma~\ref{lem:uniformBound}, i.e. $\left|\left(\operS^{-1}\right)_{n,1}\right|\leq\seriesAconst\left(\frac{l\lambda}{2}\right)^{-1}$.

The final uniform bound for this case is
\begin{equation*}
  \sum_{k=1}^N{(l+k)\left|\left(\operS^{-1}\right)_{k,1}\right|}\leq C_2\left(\frac{\lambda}{2}\right)^{-1+1/m}.
\end{equation*}

\paragraph{Case II} For $l>\alpha$.

For this case we have $d^l_1>2^{2m+1}\frac{l\lambda}{2}$, and therefore from the proof of Lemma~\ref{lem:columnEstimate} it follows that the
whole column is within the geometric decay regime, therefore the whole column can be estimated like in \eqref{eq:geometricRegime}
\begin{equation*}
  \sum_{k=1}^N{(l+k)\left|\left(\operS^{-1}\right)_{k,1}\right|\leq\left|\left(\operS^{-1}\right)_{1,1}\right|2\left(l\sum_{k=1}^\infty{\frac{1}{2^k}+\sum_{k=1}^\infty{\frac{k}{2^k}}}\right)\leq c_3\left|\left(\operS^{-1}\right)_{1,1}\right|l\leq C_3\left(\frac{\lambda}{2}\right)^{-1}}.
\end{equation*}

\paragraph{Final bound} The bound in Case I is clearly of higher order, hence it is the final uniform bound. The bound is true
for other than the first columns, as there are at most $n=[c\left(\frac{l\lambda}{2}\right)^{1/2m}]$ elements beyond the geometric decay regime.

\qed

\nocite{*}
\bibliography{CM_elliptic_system}
\bibliographystyle{alpha}

\end{document}